\documentclass[10pt]{article}
\usepackage{amsmath}
\usepackage{amsfonts}
\usepackage{amscd}
\usepackage{amsthm}
\usepackage{a4}
\usepackage{graphicx, color, epsfig}
\usepackage[english]{babel}
\usepackage{titlesec}
\titleformat{\section}{\large\bfseries}{\thesection}{1em}{}
\newtheorem{theorem}{Theorem}[section]

\newtheorem{prop}[theorem]{Proposition}

\newenvironment{PROOF}[1][Proof.]{\begin{trivlist}
\item[\hskip \labelsep {\bfseries #1}]}{\end{trivlist}}

\newcommand{\QED}{\hfill $\Box$}

\begin{document}

\numberwithin{figure}{section}
\numberwithin{equation}{section} 

\rmfamily
\setlength{\parindent}{15pt}
\title{The Kadomtsev-Petviashvili I Equation on the Half-Plane}
\author{
{D. Mantzavinos \& A. S. Fokas}\\
\\
{Department of Applied Mathematics and Theoretical Physics,}\\
{University of Cambridge, Cambridge CB3 0WA, UK.}}

\maketitle

\begin{abstract}
A new method for the solution of initial-boundary value problems for \textit{linear} and \textit{integrable nonlinear} evolution PDEs in one spatial dimension was introduced by one of the authors in 1997 \cite{F1997}. This approach was subsequently extended to initial-boundary value problems for evolution PDEs in two spatial dimensions, first in the case of linear PDEs \cite{F2002b} and, more recently, in the case of integrable nonlinear PDEs, for the Davey-Stewartson and the Kadomtsev-Petviashvili II equations on the half-plane (see \cite{FDS2009} and \cite{MF2011} respectively). In this work, we study the analogous problem for the Kadomtsev-Petviashvili I equation; in particular, through the simultaneous spectral analysis of the associated Lax pair via a d-bar formalism, we are able to obtain an integral representation for the solution, which involves certain transforms of all the initial and the boundary values, as well as an identity, the so-called global relation, which relates these transforms in appropriate regions of the complex spectral plane.
\end{abstract}

\section{Introduction}
The Kadomtsev-Petviashvili equations are two integrable nonlinear evolution PDEs in two spatial dimensions,
\begin{equation}\label{kpintro}
q_t+6qq_x+q_{xxx}+3 \sigma\partial_x^{-1} q_{yy}=0,\quad \sigma=\pm 1,
\end{equation}
where the nonlocal operator $\partial^{-1}_x$ is defined as
\begin{equation}\label{-1intro}
\partial_x^{-1}f(x)=\int_{-\infty}^{x} f(\xi)\ d\xi. 
\end{equation}
The case $\sigma=1$ is known as the KPII equation and, in the context of fluid mechanics, appears in the study of long waves in shallow water under weak surface tension. The case $\sigma=-1$ is called the KPI equation and can be employed to model water waves in thin films, where the very high surface tension dominates the gravitational force. More generally, equation \eqref{kpintro} is considered as one of the most generic models in nonlinear wave theory, modeling waves in acoustics, ferromagnetics and, more recently, in the Bose-Einstein condensation and in string theory.

The KP equations are regarded as the natural extension of the celebrated Korteweg-de Vries (KdV) equation from one to two spatial dimensions; in fact, they were first introduced in a 1970 paper by B. Kadomtsev and V. Petviashvili in order to study the stability of the soliton solutions of the KdV under the effect of transverse perturbations. KPI and KPII differ significantly not only in terms of the physical phenomena that they describe, but also in their underlying mathematical structure.

The \textit{initial value} problem for the KdV equation was solved via the Inverse Scattering transform (IST) in 1967 \cite{GGKM}; the initial value problems for the KPI and the KPII equations were formally solved in \cite{FAbl1983} and \cite{AblBarF1983}, using a nonlocal Riemann-Hilbert formalism and a d-bar formalism respectively (see also \cite{BLP1989}-\cite{Fkp2009}). A unified transform method for solving \textit{initial-boundary value} problems for linear and integrable nonlinear evolution PDEs in one spatial dimension was introduced by one of the authors in \cite{F1997}. In particular, interesting results for the linearised version of the KdV and for the KdV itself formulated on the half-line are presented in \cite{F2002a}-\cite{FTreharne2008}. The generalisation of these results from one to two spatial dimensions for linear and for integrable nonlinear equations is presented in \cite{F2002b} and \cite{FDS2009}, respectively.

An analytical approach to the initial-boundary value problem for the KPII equation on the half-plane $\{-\infty<x<\infty,\ 0<y<\infty\}$ was presented in \cite{MF2011}. Here, we will employ the general methodology of \cite{FDS2009} and \cite{MF2011} in order to analyse the KPI equation for $y$ on the half-line, that is
\begin{equation}\label{kpiintro}
q_t+6qq_x+q_{xxx}-3\partial_x^{-1} q_{yy}=0, \quad -\infty<x<\infty, \ 0<y<\infty, \ t>0.
\end{equation}
The analysis, which is proposed in sections 2-4, involves the following steps:
\begin{enumerate}
\item \textit{The formulation of the PDE in terms of a Lax pair.} For the spectral variable $k=k_R+i k_I,\ k_R,\, k_I \in \mathbb R$, and for some scalar function $\mu(x,y,t,k_R,k_I)$, the KPI equation is the compatibility condition of the Lax pair:
\begin{subequations}\label{laxpairkpi}
\begin{equation}
\mu_y+i \mu_{xx}-2k\mu_{x}=-i q \mu \label{lax1}
\end{equation}
\begin{equation}
\mu_t+4\mu_{xxx}+12ik\mu_{xx}-12k^2\mu_x=F\mu\label{lax2},
\end{equation}
\end{subequations}
where $\mu$ is bounded for all $k \in \mathbb C$ and such that
\begin{equation}\label{muasy}
\lim_{|k|\rightarrow \infty}\mu=1+\mathcal O\left(\frac{1}{k}\right)
\end{equation}
with the operator $F$ defined by
\begin{equation}\label{Fi}
F(x,y,t,k)=-6q\left(ik+\partial_x\right)-3\left(q_x+i \partial_x^{-1}q_y\right).
\end{equation}

\item \textit{The direct problem.} By applying a Fourier transform in $x$, it is possible to analyse the two equations defining the Lax pair \textit{simultaneously} in order to obtain an expression for $\mu$ which is bounded for all $k \in \mathbb C$. It turns out that $\mu$ has different representations in different parts of the complex $k$-plane, namely

\begin{figure}[ht]
\begin{center}
\resizebox{4.5cm}{!}{\input{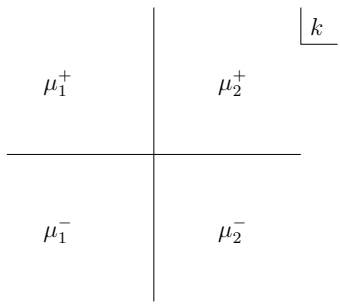}}
\end{center}
\caption{The eigenfunction $\mu$ in the four quadrants of the complex $k$-plane.}
\label{muintro}
\end{figure} 

\begin{align}\label{spectralfunctionsintro}
\mu(x,y,t,k_R,k_I) = \left\lbrace \begin{array}{ll} \mu_{2}^{+}(x,y,t,k_R,k_I), &\quad k\in I, \\ \\ \mu_{1}^{+}(x,y,t,k_R,k_I), &\quad k\in II,  \\ \\ \mu_{1}^{-}(x,y,t,k_R,k_I), &\quad k\in III, \\ \\ \mu_{2}^{-}(x,y,t,k_R,k_I), &\quad k\in IV, \end{array} \right.
\end{align}
where $I-IV$ denote the four quadrants of the complex $k$-plane.

The functions $\mu_{1,2}^\pm$ satisfy linear integral equations of Volterra type, which depend on $q(x,y,t)$, $q(x,0,t)$ and $q_y(x,0,t)$.

\item \textit{The derivation of the global relation.} This relation is an algebraic equation coupling the so-called spectral functions, which are appropriate nonlinear transforms of the initial condition and the boundary values:
\begin{equation}\label{ibviintro}
 q_0(x,y)=q(x,y,0),\ g(x,t)=q(x,0,t),\ h(x,t)=q_y(x,0,t).
\end{equation} 
The spectral functions are defined via certain linear integral equations which depend on the functions in \eqref{ibviintro}.

\item \textit{The inverse problem.} By using the fact that the function $\mu$ defined in step 2 is bounded for all $k\in \mathbb C$, it is possible to obtain an alternative representation for this function using a d-bar formalism (or more precisely the so-called Pompeiu's formula). In order to achieve this, it is necessary to: (a) compute $\partial \mu/\partial \bar k$, (b) compute the jumps of $\mu_{1,2}^\pm$ across the real and imaginary $k$-axes (actually  the jump across the imaginary $k$-axis vanishes). The d-bar derivatives and the above jumps can be expressed in terms of the spectral functions. Thus, $\mu$ can be expressed via the spectral functions and hence via $q_0$, $g$ and $h$. After obtaining $\mu$, it is straightforward to obtain a formula for $q$.
\end{enumerate}

Finally, in section 5 we will obtain an expression for the solution to the linearised version of the above problem in two different ways; first, by taking the linear limit in the representation obtained in the nonlinear case and by applying the method of \cite{F1997} for linear evolution PDEs on the linearised KP equation:
\begin{equation}\label{lkpiintro}
u_t+u_{xxx}-3\partial_x^{-1} u_{yy}=0.
\end{equation}

\paragraph{Notations and Assumptions}

\begin{itemize}
\item The complex variable $k$ is defined as
\begin{equation*}
 k = k_R + ik_I, \quad k_R,\, k_I \in \mathbb R;
\end{equation*}
the variable $l$ is real, $l \in \mathbb R$.
\item A bar on top of a complex variable will denote the complex conjugate of this variable; in particular, $\bar k = k_R-ik_I$.
\item For the solution of the inverse problem we will make use of the so-called Pompeiu's formula: if $f(x,y)$ is a smooth function in some domain $ \mathcal D\subset \mathbb R^2$ with a piece-wise smooth boundary $\partial \mathcal D$, then $f$ is related to its value on $\partial \mathcal D$ and to its d-bar derivative inside $\mathcal D$ via the equation
\begin{equation}\label{pompeiuiintro}
f(x,y)=\frac{1}{2i\pi}\int_{\partial \mathcal D}\frac{d\zeta}{\zeta-z}\, f(\xi,\eta)+\frac{1}{2i\pi}\iint_{\mathcal D}\frac{d\zeta\wedge d\bar\zeta}{\zeta-z}\, \frac{\partial f}{\partial \overline \zeta}\,(\xi, \eta), \quad \zeta=\xi+i \eta,
\end{equation}
where 
\begin{equation*}
 d \zeta \wedge d \bar \zeta=-2i d \xi d \eta.
\end{equation*}
\item We will denote the initial value by
\begin{subequations}
\begin{equation}
q(x,y,0)=q_0(x,y) \label{icintro}
\end{equation}
and the boundary values $q(x,0,t)$ and $q_y(x,0,t)$ by
\begin{equation}
q(x,0,t)=g(x,t)\label{lkpbc1intro}
\end{equation}
and
\begin{equation}\label{lkpbc2intro}
q_y(x,0,t)=h(x,t).
\end{equation}
\end{subequations}
We will assume that $q_0\in \mathbb S(\mathbb R \times \mathbb R^+)$, where $\mathbb S$ denotes the space of Schwartz functions.
\item We will seek a solution which decays as $y\rightarrow \infty$ for all fixed $(x,t)$ and which also decays as $|x|\rightarrow \infty$ for all fixed $(y,t)$.
\item Throughout this paper we will {\it assume} that there exists a solution $q(x,y,t)$ with sufficient smoothness and decay in $\bar \Omega$, which denotes the closure of the domain $\Omega$.
\item A hat ``$\wedge$''above a function will denote the Fourier transform of this function in the variable $x$.
\item The index ``o'' on a function will denote evaluation at $t=0$.
\item Whenever we write $q$, we mean that $q$ depends on the physical variables $(x,y,t)$.
\item Whenever we write $\mu$, we mean that $\mu$ depends on the physical variables $(x,y,t)$ and on the spectral variables $k_R$, $k_I\in \mathbb R$.
\item The subscripts $X$ and $t$ on a function $f$ denote its dependence on $(x,y,t)$, i.e. 
\begin{equation*} 
 f_{X t}=f(x,y,t).
\end{equation*}
Similarly, 
\begin{equation*} 
 f_{\Xi t}=f(\xi,\eta,t).
\end{equation*}
\end{itemize}

\section{The direct problem}

\begin{prop}\label{kpidpprop}
Assume that there exists a solution $q(x,y,t)$,  defined for $(x,y,t)\in \Omega$, to a well-posed initial-boundary value problem for equation \eqref{kpiintro}. Then, the function $\mu$, which is bounded for all $k \in \mathbb C$ and is defined as
\begin{align}\label{spectralfunctionsi}
\mu(x,y,t,k_R,k_I) = \left\lbrace \begin{array}{ll} \mu_{1}^{+}(x,y,t,k_R,k_I), &\quad k_R\leq0,\ k_I\geq0, \\ \\ \mu_{1}^{-}(x,y,t,k_R,k_I), &\quad k_R\leq0,\ k_I\leq0,  \\ \\ \mu_{2}^{-}(x,y,t,k_R,k_I), &\quad k_R\geq0,\ k_I\leq0, \\ \\ \mu_{2}^{+}(x,y,t,k_R,k_I), &\quad k_R\geq0,\ k_I\geq0, \end{array} \right.
\end{align}
admits the following representations in terms of $q(x,y,t)$, $q(x,0,t)$ and $q_y(x,0,t)$:
\begin{subequations}\label{plusminus12i}
\begin{align}\label{plus1i}
\mu_1^+&=1-\frac{i}{2\pi}\int_{0}^{\infty}dl\int_{-\infty}^{\infty}d\xi\int_{0}^{y}d\eta\, e^{-il(\xi-x)-il(l+2k)(\eta-y)}q\mu_1^+\nonumber\\
&-\frac{1}{2\pi}\int_{0}^{-2k_R}dl\int_{-\infty}^{\infty}d\xi\int_{t}^Td\tau\,e^{-il(\xi-x)+il(l+2k)y+\omega(k,l)(\tau-t)}H\phi_1^+\nonumber\\
&+\frac{1}{2\pi}\int_{-2k_R}^{\infty}dl\int_{-\infty}^{\infty}d\xi\int_{0}^{t}d\tau\,e^{-il(\xi-x)+il(l+2k)y+\omega(k,l)(\tau-t)}H\phi_1^+\nonumber\\
&+\frac{i}{2\pi}\int_{-\infty}^{0}dl\int_{-\infty}^{\infty}d\xi\int_{y}^{\infty}d\eta\,e^{-il(\xi-x)-il(l+2k)(\eta-y)}q\mu_1^+, \quad k_R\leq 0,\, k_I\geq 0,
\end{align}
\begin{align}\label{minus1i}
\mu_1^-&=1+\frac{i}{2\pi}\int_{0}^{\infty}dl\int_{-\infty}^{\infty}d\xi\int_{y}^\infty d\eta\, e^{-il(\xi-x)-il(l+2k)(\eta-y)}q\mu_1^-\nonumber\\
&-\frac{1}{2\pi}\int_{-\infty}^0dl\int_{-\infty}^{\infty}d\xi\int_{t}^{T}d\tau\,e^{-il(\xi-x)+il(l+2k)y+\omega(k,l)(\tau-t)}H\phi_1^-\nonumber\\
&-\frac{i}{2\pi}\int_{-\infty}^{0}dl\int_{-\infty}^{\infty}d\xi\int_{0}^{y}d\eta\,e^{-il(\xi-x)-il(l+2k)(\eta-y)}q\mu_1^-,\quad  k_R\leq 0,\, k_I\leq 0,
\end{align}
\begin{align}\label{minus2i}
\mu_2^-&=1+\frac{i}{2\pi}\int_{0}^{\infty}dl\int_{-\infty}^{\infty}d\xi\int_{y}^\infty d\eta\, e^{-il(\xi-x)-il(l+2k)(\eta-y)}q\mu_2^-\nonumber\\
&-\frac{1}{2\pi}\int_{-\infty}^0dl\int_{-\infty}^{\infty}d\xi\int_{t}^{T}d\tau\,e^{-il(\xi-x)+il(l+2k)y+\omega(k,l)(\tau-t)}H\phi_2^-\nonumber\\
&+\frac{1}{2\pi}\int_{-2k_R}^0dl\int_{-\infty}^{\infty}d\xi\int_{0}^{t}d\tau\,e^{-il(\xi-x)+il(l+2k)y+\omega(k,l)(\tau-t)}H\phi_2^-\nonumber\\
&-\frac{i}{2\pi}\int_{-\infty}^{0}dl\int_{-\infty}^{\infty}d\xi\int_{0}^{y}d\eta\,e^{-il(\xi-x)-il(l+2k)(\eta-y)}q\mu_2^-,\quad  k_R\geq 0,\, k_I\leq 0,
\end{align}
and
\begin{align}\label{plus2i}
\mu_2^+&=1-\frac{i}{2\pi}\int_{0}^{\infty}dl\int_{-\infty}^{\infty}d\xi\int_{0}^{y}d\eta\, e^{-il(\xi-x)-il(l+2k)(\eta-y)}q\mu_2^+\nonumber\\
&+\frac{1}{2\pi}\int_{0}^{\infty}dl\int_{-\infty}^{\infty}d\xi\int_{0}^{t}d\tau\,e^{-il(\xi-x)+il(l+2k)y+\omega(k,l)(\tau-t)}H\phi_2^+\nonumber\\
&+\frac{i}{2\pi}\int_{-\infty}^{0}dl\int_{-\infty}^{\infty}d\xi\int_{y}^{\infty}d\eta\,e^{-il(\xi-x)-il(l+2k)(\eta-y)}q\mu_2^+, \quad k_R\geq 0,\, k_I\geq 0.
\end{align}
\end{subequations}
where
\begin{align}
\omega(k,l)&:=-4il(l^2+3kl+3k^2)\label{omega},\\
\phi^\pm_{j}(x,t,k_R,k_I)&:=\mu^\pm_{j}(x,0,t,k_R,k_I),\quad j=1,2,\label{phidefi}\\
H(x,t,k,l)&:=3\left[g_x(x,t)-2i(l+k)g(x,t)-i\partial_x^{-1}h(x,t)\right]\label{H},
\end{align}
and $g, h$ are defined in \eqref{ibviintro}. 

\end{prop}

\begin{PROOF}  

Define the Fourier transform of $\mu$ with respect to $x$ by
\begin{equation}\label{ft}
\hat \mu(l,y,t,k_R,k_I)=\int_{-\infty}^{\infty}dx\, e^{-ilx}\Big[\mu(x,y,t,k_R,k_I)-1\Big],
\end{equation}
with the inverse given by
\begin{equation}\label{ift}
\mu(x,y,t,k_R,k_I)=1+\frac{1}{2\pi}\int_{-\infty}^{\infty}dl\, e^{ilx}\hat\mu(l,y,t,k_R,k_I).
\end{equation}
Applying \eqref{ft} to the first equation of the Lax pair yields:
\begin{equation}\label{intfac}
\left(\hat \mu e^{-il(l+2k)y}\right)_y=-i e^{-il(l+2k)y}\int_{-\infty}^{\infty}dx\, e^{-ilx}{q\mu}.
\end{equation}
This equation can be integrated either from 0 to $y$ or from $y$ to $\infty$, depending on the sign of the real part of the exponent $-il(l+2k)$, namely
\begin{align}\label{intfac2}
\hat{\mu}(l,y,t,k_R,k_I)=\left\{\begin{array}{ll} \!\!\!\!-i\int_0^y d\eta\, \int_{-\infty}^{\infty}d\xi\, e^{-il\xi-il(l+2k)(\eta-y)}{q\mu}+\hat{\mu}(l,0,t,k_R,k_I)\, e^{il(l+2k)y},\ lk_I \geq 0, \\ \\ i\int_y^\infty d\eta\int_{-\infty}^{\infty}d\xi\, e^{-il\xi-il(l+2k)(\eta-y)}{q\mu}, \ lk_I \leq 0. \end{array}\right.
\end{align}
Hence, equation \eqref{ift} implies the following formulae for $\mu^+$ and $\mu^-$, defined for $k_I\geq 0$ and $k_I\leq 0$ respectively:
\begin{subequations}\label{plusminusi}
{\small\begin{align}\label{plusi}
&\mu^+(x,y,t,k_R,k_I)=1-\frac{i}{2\pi}\int_{0}^{\infty}dl\int_{-\infty}^{\infty}d\xi\int_{0}^{y}d\eta\, e^{-il(\xi-x)-il(l+2k)(\eta-y)}\mbox{\large $($}q\mu^+\mbox{\large $)$}(\xi,\eta,t,k_R,k_I)\nonumber\\
&+\frac{1}{2\pi}\int_{0}^{\infty}dl\,e^{ilx+il(l+2k)y}\hat\phi^+(l,t,k_R,k_I)\nonumber\\
&+\frac{i}{2\pi}\int_{-\infty}^{0}dl\int_{-\infty}^{\infty}d\xi\int_{y}^{\infty}d\eta\,e^{-il(\xi-x)-il(l+2k)(\eta-y)}\mbox{\large $($}q\mu^+\mbox{\large $)$}(\xi,\eta,t,k_R,k_I), \quad k_I\geq 0,
\end{align}}
and
{\small\begin{align}\label{minusi}
&\mu^-(x,y,t,k_R,k_I)=1+\frac{i}{2\pi}\int_{0}^{\infty}dl\int_{-\infty}^{\infty}d\xi\int_{y}^\infty d\eta\, e^{-il(\xi-x)-il(l+2k)(\eta-y)}\mbox{\large $($}q\mu^-\mbox{\large $)$}(\xi,\eta,t,k_R,k_I)\nonumber\\
&+\frac{1}{2\pi}\int_{-\infty}^0dl\,e^{ilx+il(l+2k)y}\hat\phi^-(l,t,k_R,k_I)\nonumber\\
&-\frac{i}{2\pi}\int_{-\infty}^{0}dl\int_{-\infty}^{\infty}d\xi\int_{0}^{y}d\eta\,e^{-il(\xi-x)-il(l+2k)(\eta-y)}\mbox{\large $($}q\mu^-\mbox{\large $)$}(\xi,\eta,t,k_R,k_I),\quad k_I\leq 0,
\end{align}}
\end{subequations}
where $\hat \phi^\pm$ are the Fourier transforms of $\phi^\pm$, defined by equations \eqref{phidefi}.

In order to compute $\hat \phi^\pm$, we use the $t$-part of the Lax pair, namely we evaluate equation \eqref{lax2} at $y=0$ and then apply the Fourier transform \eqref{ft} to the resulting equation:
\begin{equation}\label{phieq}
\hat \phi_t(l,t,k_R,k_I)+\omega(k,l)\hat \phi(l,t,k_R,k_I)=\int_{-\infty}^{\infty}d\xi\, e^{-il\xi}H(\xi,t,k,l)\phi(\xi,t,k_R,k_I), 
\end{equation}
where $\omega$ and $H$ are defined by \eqref{omega} and \eqref{H} respectively.

Equation \eqref{phihat} can be integrated either from $0$ to $t$ or from $t$ to $T$, depending on the sign of the real part of $\omega$, i.e. on the sign of $lk_I(l+2k_R)$, see equation \eqref{omega}:
\begin{align}
\hat \phi(l,t,k_R,k_I)=\left\{\begin{array}{ll}&\!\!\!\!\!\!\!\!\int_0^td\tau\int_{-\infty}^{\infty}d\xi\, e^{-il\xi+\omega(k,l)(\tau-t)}H\phi+e^{-\omega(k,l)t}\hat\phi(l,0,k_R,k_I),\ lk_I(l+2k_R)\geq 0\\ \\ &\!\!\!\!\!\!\!\!-\int_t^Td\tau\int_{-\infty}^{\infty}d\xi\, e^{-il\xi+\omega(k,l)(\tau-t)}H\phi+e^{\omega(k,l)(T-t)}\hat\phi(l,T,k_R,k_I),\ lk_I(l+2k_R)\leq 0.\end{array}\right.\label{phihat}
\end{align}
For $\hat \phi^+$ we have that $k_I\geq0$. We distinguish two subcases, namely $k_R\geq0$ and $k_R\leq0$:
\vskip 2mm
\noindent\underline{$k_I\geq 0$ and $k_R\geq 0$} 
\vskip 2mm
\noindent For $l\in(-\infty,-2k_R)\cup(0,\infty)$ we have that $l(l+2k_R)\geq 0$ thus we use the first expression in \eqref{phihat}, whereas for $l\in(-2k_R,0)$ we have that $l(l+2k_R)\leq 0$ and so we use the second expression in \eqref{phihat}:
\begin{subequations}\label{phihat2}
{\small\begin{align}
\hat \phi_2^+(l,t,k_R,k_I)=\left\{\begin{array}{ll}&\!\!\!\!\!\!\!\!\int_0^td\tau\int_{-\infty}^{\infty}d\xi\, e^{-il\xi+\omega(k,l)(\tau-t)}H\phi_2^++e^{-\omega(k,l)t}\hat\phi_2^+(l,0,k_R,k_I),\ l\in(-\infty,-2k_R]\cup[0,\infty)\\ \\ &\!\!\!\!\!\!\!\!-\int_t^Td\tau\int_{-\infty}^{\infty}d\xi\, e^{-il\xi+\omega(k,l)(\tau-t)}H\phi_2^++e^{\omega(k,l)(T-t)}\hat\phi_2^+(l,T,k_R,k_I),\ l\in[-2k_R,0],\end{array}\right.\label{phihat2+}
\end{align}}
and
\vskip 2mm
\noindent\underline{$k_I\geq 0$ and $k_R\leq 0$} 
\vskip 2mm
\noindent For $l\in(-\infty,0)\cup(-2k_R,\infty)$ we have that $l(l+2k_R)\geq 0$ thus we use the first expression in \eqref{phihat}, whereas for $l\in(0,-2k_R)$ we have that $l(l+2k_R)\leq 0$ and so we use the second expression in \eqref{phihat}:
{\small\begin{align}
\hat \phi_1^+(l,t,k_R,k_I)=\left\{\begin{array}{ll}&\!\!\!\!\!\!\!\!\int_0^td\tau\int_{-\infty}^{\infty}d\xi\, e^{-il\xi+\omega(k,l)(\tau-t)}H\phi_1^++e^{-\omega(k,l)t}\hat\phi_1^+(l,0,k_R,k_I),\ l\in(-\infty,0]\cup[-2k_R,\infty)\\ \\ &\!\!\!\!\!\!\!\!-\int_t^Td\tau\int_{-\infty}^{\infty}d\xi\, e^{-il\xi+\omega(k,l)(\tau-t)}H\phi_1^++e^{\omega(k,l)(T-t)}\hat\phi_1^+(l,T,k_R,k_I),\ l\in[0,-2k_R].\end{array}\right.\label{phihat1+}
\end{align}}
In the case of $\phi^-$, we have that $k_I\leq0$. As before, we distinguish two subcases, $k_R\leq0$ and $k_R\geq0$:
\vskip 2mm
\noindent\underline{$k_I\leq 0$ and $k_R\leq 0$} 
\vskip 2mm
\noindent For $l\in(-\infty,0)\cup(-2k_R,\infty)$ we have that $l(l+2k_R)\geq 0$ thus we use the second expression in \eqref{phihat}, whereas for $l\in(0,-2k_R)$ we have that $l(l+2k_R)\leq 0$ and so we use the first expression in \eqref{phihat}:
{\small\begin{align}
\hat \phi_1^-(l,t,k_R,k_I)=\left\{\begin{array}{ll}&\!\!\!\!\!\!\!\!-\int_t^Td\tau\int_{-\infty}^{\infty}d\xi\, e^{-il\xi+\omega(k,l)(\tau-t)}H\phi_1^-+e^{\omega(k,l)(T-t)}\hat\phi_1^-(l,T,k_R,k_I),\ l\in(-\infty,0]\cup[-2k_R,\infty) \\ \\&\!\!\!\!\!\!\!\!\int_0^td\tau\int_{-\infty}^{\infty}d\xi\, e^{-il\xi+\omega(k,l)(\tau-t)}H\phi_1^-+e^{-\omega(k,l)t}\hat\phi_1^+(l,0,k_R,k_I),\ l\in[0,-2k_R], \end{array}\right.\label{phihat1-}
\end{align}}
and
\vskip 2mm
\noindent\underline{$k_I\leq 0$ and $k_R\geq 0$} 
\vskip 2mm
\noindent For $l\in(-\infty,-2k_R)\cup(0,\infty)$ we have that $l(l+2k_R)\geq 0$ thus we use the second expression in \eqref{phihat}, whereas for $l\in(-2k_R,0)$ we have that $l(l+2k_R)\leq 0$ and so we use the first expression in \eqref{phihat}:
{\small\begin{align}
\hat \phi_2^-(l,t,k_R,k_I)=\left\{\begin{array}{ll}&\!\!\!\!\!\!\!\!-\int_t^Td\tau\int_{-\infty}^{\infty}d\xi\, e^{-il\xi+\omega(k,l)(\tau-t)}H\phi_2^-+e^{\omega(k,l)(T-t)}\hat\phi_2^-(l,T,k_R,k_I),\ l\in(-\infty,-2k_R]\cup[0,\infty) \\ \\ &\!\!\!\!\!\!\!\!\int_0^td\tau\int_{-\infty}^{\infty}d\xi\, e^{-il\xi+\omega(k,l)(\tau-t)}H\phi_2^-+e^{-\omega(k,l)t}\hat\phi_2^-(l,0,k_R,k_I),\ l\in[-2k_R,0]. \end{array}\right.\label{phihat2-}
\end{align}}
\end{subequations}
Moreover, equations \eqref{plusminusi} evaluated at $y=0$ imply that $\phi^\pm$ satisfy the following equations: 
\begin{subequations}\label{choice}
\begin{equation}\label{phi+choice}
\phi^+=1+\frac{1}{2\pi}\int_0^\infty\!\! dl\, e^{ilx}\hat \phi^++\frac{i}{2\pi}\int_{-\infty}^{0}\!\!dl\int_{-\infty}^{\infty}\!\!d\xi\int_0^\infty \!\!d\eta\, e^{-il(\xi-x)-il(l+2k)\eta}q\mu^+,
\end{equation}
\begin{equation}\label{phi-choice}
\phi^-=1+\frac{i}{2\pi}\int_{0}^{\infty}\!\!dl\int_{-\infty}^{\infty}\!\!d\xi\int_0^\infty \!\!d\eta\, e^{-il(\xi-x)-il(l+2k)\eta}q\mu^-+\frac{1}{2\pi}\int_{-\infty}^0\!\! dl\, e^{ilx}\hat \phi^-.
\end{equation}
\end{subequations}
Hence, together with the inverse Fourier transform \eqref{ift} for $\phi^\pm$, the above imply that $\hat \phi^+$ has no restrictions on $t$ for $l\geq 0$ as well as that $\hat \phi^-$ has no restrictions on $t$ for $l\leq 0$. Therefore we can choose the values of $\phi^\pm$ for $t=0$ and $t=T$ appropriately as:
\begin{subequations}\label{choice2}
\begin{equation}
\hat\phi_1^+(l,T,k_R,k_I)=0, \quad l\in[0,-2k_R], \quad \hat\phi_1^+(l,0,k_R,k_I)=0, \quad l\in[-2k_R,\infty),
\end{equation}
\begin{equation}
\hat\phi_1^-(l,T,k_R,k_I)=0, \quad l\in(-\infty,0],
\end{equation}
\begin{equation}
\hat\phi_2^-(l,T,k_R,k_I)=0, \quad\in(-\infty,-2k_R], \quad \hat\phi_2^-(l,0,k_R,k_I)=0, \quad l\in[-2k_R,0],
\end{equation}
and
\begin{equation}
\hat\phi_2^+(l,0,k_R,k_I)=0, \quad l\in[0,\infty). 
\end{equation}
\end{subequations}
Using these equations in equations \eqref{phihat2} and substituting the resulting expressions in \eqref{plusminusi}, we find equations $\eqref{plus1i}-\eqref{plus2i}$. \QED
\end{PROOF}

\section{The global relation}

\begin{prop}\label{griprop}
Define $q_0$, $g$ and $h$ by equation \eqref{ibviintro}. Also, define the functions $\omega$, $\phi$ and $H$ by equations \eqref{omega}, \eqref{phidefi} and \eqref{H} respectively. Then, the following identity, called the \textit{global relation}, is valid for $\mu$:
\begin{align}
&\int_{-\infty}^{\infty}d\xi\int_{0}^\infty d\eta\, e^{-il\xi-il(l+2k)\eta} q_0(\xi,\eta)\mu(\xi,\eta,0,k_R,k_I)\nonumber\\
&-i\int_{-\infty}^{\infty}d\xi\int_{0}^{t}d\tau\, e^{-il\xi+\omega(k,l)\tau}H(\xi,\tau,k,l)\phi(\xi,\tau,k_R,k_I)\nonumber\\
&=\int_{-\infty}^{\infty}d\xi\int_{0}^\infty d\eta\, e^{-il\xi-il(l+2k)\eta+\omega(k,l)t} q(\xi,\eta,t)\mu(\xi,\eta,t,k_R,k_I),\quad k \in \mathbb C,\quad lk_I\leq0. \label{gri}
\end{align}
Moreover, an alternative form of the global relation is:
\begin{align}
&\int_{-\infty}^{\infty}d\xi\int_{0}^\infty d\eta\, e^{-il\xi-il(l+2k)\eta} q(\xi,\eta,t)\mu(\xi,\eta,t,k_R,k_I)\nonumber\\
&-i\int_{-\infty}^{\infty}d\xi\int_{t}^{T}d\tau\, e^{-il\xi+\omega(k,l)(\tau-t)}H(\xi,\tau,k,l)\phi(\xi,\tau,k_R,k_I)\nonumber\\
&=\int_{-\infty}^{\infty}d\xi\int_{0}^\infty d\eta\, e^{-il\xi-il(l+2k)\eta+\omega(k,l)(T-t)} q(\xi,\eta,T)\mu(\xi,\eta,T,k_R,k_I),\quad k \in \mathbb C,\quad lk_I\leq0. \label{gr2i}
\end{align}
\end{prop}

\begin{PROOF}

Applying the Fourier transform \eqref{ft} to the second Lax equation yields an ODE in $t$,
\begin{equation}\label{mut}
\left(\hat \mu\, e^{\omega(k,l)t}\right)_t=e^{\omega(k,l)t}\int_{-\infty}^{\infty}d\xi\, e^{-il\xi}F\mu, 
\end{equation}
where $F$ is defined by \eqref{Fi}.

This equation, together with equation \eqref{intfac}, imply the following equality:
\begin{equation}\label{derivatives}
\!\!\left(\!e^{-il(l+2k)y+\omega(k,l)t}\int_{-\infty}^{\infty}d\xi\, e^{-il\xi}q\mu\!\right)_t\!\!=\!\left(\!i\, e^{-il(l+2k)y+\omega(k,l)t}\!\!\int_{-\infty}^{\infty}\!\!d\xi\, e^{-il\xi}F\mu\!\right)_y.
\end{equation}
Then, Green's theorem in the plane, 
\begin{equation}
\iint_{\mathcal D}\left(\frac{\partial g}{\partial u}-\frac{\partial f}{\partial v}\right)du dv=\int_{\partial \mathcal D}f(u,v)du+g(u,v)dv,
\end{equation}
applied inside the domains depicted in figures \ref{green1i} and \ref{green2i}, yields the identities \eqref{gri} and \eqref{gr2i}, respectively.
\begin{figure}[ht]
\begin{center}
\resizebox{5cm}{!}{\input{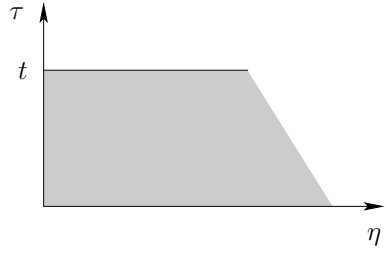}}
\end{center}
\caption{The domain of integration for Green's theorem.}
\label{green1i}
\end{figure} 

\begin{figure}[ht]
\begin{center}
\resizebox{5cm}{!}{\input{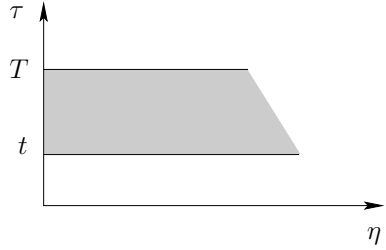}}
\end{center}
\caption{An alternative domain of integration for Green's theorem.}
\label{green2i}
\end{figure} 
Note that the constraint on $l$ is needed in order for the exponential $e^{-il(l+2k)\eta}$ to be bounded when $\eta \rightarrow \infty$.
\QED
\end{PROOF}

\section{The inverse problem}

\begin{figure}[ht]
\begin{center}
\resizebox{5cm}{!}{\input{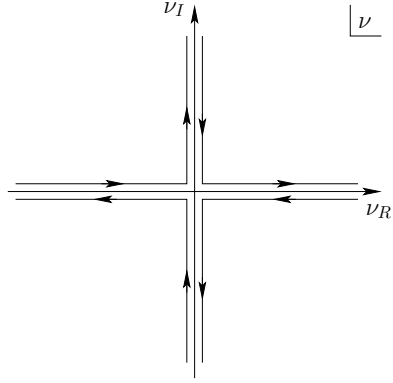}}
\end{center}
\caption{The contour of integration in Pompeiu's formula.}
\label{pompifig}
\end{figure}

Define the functions
\begin{equation}\label{deltamu}
\delta \mu^+(x,y,t,k_I)=\left(\mu_1^+-\mu_2^+\right)\!\Big|_{k_R=0}\ , \quad \delta \mu^-(x,y,t,k_I)=\left(\mu_1^--\mu_2^-\right)\!\Big|_{k_R=0}
\end{equation}
and
\begin{equation}\label{Deltamu}
\Delta \mu_1(x,y,t,k_R)=\left(\mu_1^+-\mu_1^-\right)\!\Big|_{k_I=0}\ , \quad \Delta \mu_2(x,y,t,k_R)=\left(\mu_2^+-\mu_2^-\right)\!\Big|_{k_I=0}.
\end{equation}
Employing Pompeiu's formula \eqref{pompeiuiintro} in the domain depicted in figure \ref{pompifig},

we can express the function $\mu$ in the following form:
\begin{align}\label{pompi}
\mu&=1+\frac{1}{2i \pi} \int_{0}^{\infty} \frac{d\nu_R}{\nu_R-k}\ \Delta\mu_2+\frac{1}{2i \pi} \int_{-\infty}^0 \frac{d\nu_R}{\nu_R-k}\ \Delta \mu_1\nonumber\\
&+\frac{1}{2 \pi} \int_{0}^{\infty} \frac{d\nu_I}{i \nu_I-k}\ \delta\mu^+ +\frac{1}{2 \pi} \int_{-\infty}^0 d\nu_I\ \frac{\delta \mu^-}{i\nu_I-k}\nonumber\\
&-\frac{1}{\pi}\int_{-\infty}^{0}d\nu_R\int_0^\infty \frac{d\nu_I}{\nu-k}\frac{\partial \mu_1^+}{\partial \bar \nu }-\frac{1}{\pi}\int_{-\infty}^{0}d\nu_R\int_{-\infty}^0 \frac{d\nu_I}{\nu-k}\frac{\partial \mu_1^-}{\partial \bar \nu }\nonumber\\
&-\frac{1}{\pi}\int_0^{\infty}d\nu_R\int_{-\infty}^0 \frac{d\nu_I}{\nu-k}\frac{\partial \mu_2^-}{\partial \bar \nu }-\frac{1}{\pi}\int_0^{\infty}d\nu_R\int_0^\infty \frac{d\nu_I}{\nu-k}\frac{\partial \mu_2^+}{\partial \bar \nu }.
\end{align}
Therefore, we need to compute the jumps of $\mu$ across the real and the imaginary $k$-axes as well as the d-bar derivatives in the four quandrants of the complex $k$-plane.

\begin{prop}(The d-bar derivatives)\label{kpid-barprop}
Define the functions $\mu_j^\pm, \ j=1,2,$ by equations \eqref{spectralfunctionsi} and $\eqref{plus1i}-\eqref{plus2i}$. Then,
\begin{subequations}
\begin{equation}\label{dbarplus1i}
\frac{\partial \mu_1^+}{\partial \bar k}(x,y,t,k_R,k_I)=e_{1_{Xt}}f_1^+(k_R,k_I)\mu_2^+(x,y,t,-k_R,k_I),\quad k_R\leq0,\, k_I\geq 0,
\end{equation}
and
\begin{equation}\label{dbarminus2i}
\frac{\partial \mu_2^-}{\partial \bar k}(x,y,t,k_R,k_I)=e_{1_{Xt}}f_2^-(k_R,k_I)\mu_1^-(x,y,t,-k_R,k_I), \quad k_R\geq0,\, k_I\leq 0,
\end{equation}
\end{subequations}
where 
\begin{equation}\label{eXtprop}
e_{1_{Xt}}=e_{1_{Xt}}(x,y,t,k_R,k_I)=e^{-2ik_R x+4k_Rk_I y-\omega(k,-2k_R)t},
\end{equation}
the functions $f_j^+,\ j=1,2,$ are defined by
\begin{subequations}\label{f}
\begin{equation}\label{fplus1i}
f_1^+(k_R,k_I)=\frac{1}{2\pi}\int_{-\infty}^{\infty}d\xi\int_{0}^{T}d\tau\,e^{2ik_R\xi+\omega(k,-2k_R)\tau}H(\xi,\tau,k,-2k_R)\phi_1^+,
\end{equation}
\begin{equation}\label{fminus2i}
f_2^-(k_R,k_I)=\frac{1}{2\pi}\int_{-\infty}^{\infty}d\xi\int_{0}^{T}d\tau\,e^{2ik_R\xi+\omega(k,-2k_R)\tau}H(\xi,\tau,k,-2k_R)\phi_2^-
\end{equation}
\end{subequations}
and the functions $\omega$, $\phi_{1}^+$, $\phi_{2}^-$ and $H$ are defined by equations $\eqref{omega}-\eqref{H}$.

Moreover,
\begin{subequations}
\begin{equation}\label{dbarminus1i}
 \frac{\partial \mu_1^-}{\partial \bar k}(x,y,t,k_R,k_I)=0, \quad k_R\leq0,\, k_I\leq 0,
\end{equation}
and
\begin{equation}\label{dbarplus2i}
 \frac{\partial \mu_2^+}{\partial \bar k}(x,y,t,k_R,k_I)=0, \quad k_R\geq0,\, k_I\geq 0.
\end{equation}
\end{subequations}
\end{prop}

\begin{PROOF} 
Applying the d-bar derivative with respect to the spectral variables on equation \eqref{plus1i} and recalling that $H$ depends on $k$ but \textit{not} on $\bar k$, we have:
\begin{align}\label{dplus1i}
\frac{\partial\mu_1^+}{\partial \bar k}&=-\frac{i}{2\pi}\int_{0}^{\infty}dl\int_{-\infty}^{\infty}d\xi\int_{0}^{y}d\eta\, e^{-il(\xi-x)-il(l+2k)(\eta-y)}q\frac{\partial\mu_1^+}{\partial \bar k}\nonumber\\
&-\frac{1}{2\pi}\int_{0}^{-2k_R}dl\int_{-\infty}^{\infty}d\xi\int_{t}^Td\tau\,e^{-il(\xi-x)+il(l+2k)y+\omega(k,l)(\tau-t)}H\frac{\partial\phi_1^+}{\partial \bar k}\nonumber\\
&+\frac{1}{2\pi}\int_{-2k_R}^{\infty}dl\int_{-\infty}^{\infty}d\xi\int_{0}^{t}d\tau\,e^{-il(\xi-x)+il(l+2k)y+\omega(k,l)(\tau-t)}H\frac{\partial\phi_1^+}{\partial \bar k}\nonumber\\
&+\frac{i}{2\pi}\int_{-\infty}^{0}dl\int_{-\infty}^{\infty}d\xi\int_{y}^{\infty}d\eta\,e^{-il(\xi-x)-il(l+2k)(\eta-y)}q\frac{\partial\mu_1^+}{\partial \bar k}\nonumber\\ 
&+e^{-2ik_Rx+4k_Rk_Iy-\omega(k,-2k_R)t}f_1^+(k_R,k_I,T),\quad k_R\leq 0,\, k_I\geq0,
\end{align}
where $f_1^+$ is defined by \eqref{fplus1i}.

In equation \eqref{intfac}, replace $k$ by $-\bar k$ and introduce the notation $\tilde \mu$ for
\begin{equation}\label{mutilde}
\tilde \mu(x,y,t,k_R,k_I)=\mu(x,y,t,-k_R,k_I).
\end{equation}
Then,
\begin{equation}\label{tildeintfac}
\left(\hat{\tilde \mu}e^{-il(l-2\bar k)y}\right)_y=-ie^{-il(l-2\bar k)}\int_{-\infty}^{\infty}d\xi\, e^{-il\xi}q \tilde\mu,
\end{equation}
hence
\begin{align}\label{tildeintfac2}
\hat{\tilde\mu}=\left\{\begin{array}{ll} \!\!\!\!-i\int_0^y d\eta\, \int_{-\infty}^{\infty}d\xi\, e^{-il\xi-il(l-2\bar k)(\eta-y)}q\tilde\mu+\hat{\tilde\phi}\, e^{il(l-2\bar k)y}, \quad lk_I\geq 0\\ \\ i\int_y^\infty d\eta\int_{-\infty}^{\infty}d\xi\, e^{-il\xi-il(l-2\bar k)(\eta-y)}q\tilde\mu,\quad lk_I\leq 0, \end{array}\right.
\end{align}
where $\tilde \phi$ denotes the evaluation of $\tilde \mu$ at $y=0$.

Let
\begin{equation}\label{eXtgen}
e_{Xt}=e_{Xt}(x,y,t,k,\lambda):=e^{-i(\lambda+2k)x+i\lambda(\lambda+2k)y+\omega(k,\lambda)t-8ik^3 t}.
\end{equation}
where the indices $X,t$ denote the dependence on the variables $(x,y,t)$. Furthermore, define
\begin{equation}\label{eXt}
e_{1_{Xt}}:=e_{Xt}(x,y,t,k,-2ik_I)=e^{-2ik_R x+4k_Rk_I y-\omega(k,-2k_R)t}.
\end{equation}
In the case of $k_R\leq 0$ and $k_I\geq0$, multiplying equation \eqref{tildeintfac2} by $e_{1_{Xt}}$ yields:
\begin{align}\label{tildeintfac3}
e_{1_{Xt}}\hat{\tilde\mu}_2^+=\left\{\begin{array}{ll} \!\!\!\!-i\int_0^y d\eta\, \int_{-\infty}^{\infty}d\xi\, e^{-il\xi+2ik_R(\xi-x)-il(l-2\bar k)(\eta-y)-4k_Rk_I(\eta-y)}q(e_{1_{\Xi t}}\tilde\mu_2^+)+e_{1_{Xt}}\hat{\tilde\phi}_2^+\, e^{il(l-2\bar k)y} \\ \\ i\int_y^\infty d\eta\int_{-\infty}^{\infty}d\xi\, e^{-il\xi+2ik_R(\xi-x)-il(l-2\bar k)(\eta-y)-4k_Rk_I(\eta-y)}q(e_{1_{\Xi t}}\tilde\mu_2^+), \end{array}\right.
\end{align}
where the index $\Xi$ indicates the dependence on the variables $(\xi, \eta)$ instead of $(x, y)$.

According to equation \eqref{ift},
\begin{align}\label{tildeintfac4}
e_{1_{Xt}}\tilde \mu_2^+=e_{1_{Xt}}+\frac{1}{2\pi}\int_{-\infty}^{\infty}dl\, e^{ilx}(e_{1_{Xt}}\hat{\tilde\mu}_2^+),
\end{align}
thus, using \eqref{tildeintfac3} into \eqref{tildeintfac4} we find:
{\small\begin{align}\label{tildeintfac5}
e_{1_{Xt}}\tilde \mu_2^+&=e_{1_{Xt}}-\frac{i}{2\pi}\int_{2k_R}^{\infty}dl\int_{-\infty}^{\infty}d\xi\int_0^y d\eta\, e^{-i(l-2k_R)(\xi-x)-(il(l-2\bar k)+4k_Rk_I)(\eta-y)}q(e_{1_{\Xi t}}\tilde\mu_2^+)\nonumber\\
&+\frac{1}{2\pi}\int_{2k_R}^{\infty}dl\,e^{ilx+il(l-2\bar k)y}( e_{1_{Xt}}\hat{\tilde \phi}_2^+)\nonumber\\
&+\frac{i}{2\pi}\int_{-\infty}^{2k_R}dl\int_{-\infty}^{\infty}d\xi\int_y^\infty d\eta\, e^{-i(l-2k_R)(\xi-x)-(il(l-2\bar k)+4k_Rk_I)(\eta-y)}q(e_{1_{\Xi t}}\tilde\mu_2^+).
\end{align}}
Replacing $k$ by $-\bar k$ in equation \eqref{phieq} and integrating with respect to $t$ yields:
\begin{align}
\hat {\tilde\phi}=\left\{\!\begin{array}{ll}&\!\!\!\!\!\!\!\!\int_0^td\tau\int_{-\infty}^{\infty}d\xi\, e^{-il\xi+\omega(-\bar k,l)(\tau-t)}H\tilde\phi+e^{-\omega(-\bar k,l)t}\hat{\tilde\phi}\big|_{t=0}, \quad \mathrm{Re}\,\omega(-\bar k,l) \geq 0\\ \\ &\!\!\!\!\!\!\!\!-\int_t^Td\tau\int_{-\infty}^{\infty}d\xi\, e^{-il\xi+\omega(-\bar k,l)(\tau-t)}H\tilde\phi+e^{\omega(-\bar k,l)(T-t)}\hat{\tilde\phi}\big|_{t=T},\  \mathrm{Re}\, \omega(-\bar k,l)\leq 0.\end{array}\right.\label{tildephihat}
\end{align}
Hence, for $k_R\leq0$ and $k_I\geq 0$, the following equation is valid:
{\small\begin{align}\label{tildemuplus2}
&e_{1_{Xt}}\tilde \mu_2^+=e_{1_{Xt}}-\frac{i}{2\pi}\int_{2k_R}^{\infty}dl\int_{-\infty}^{\infty}d\xi\int_0^y d\eta\, e^{-i(l-2k_R)(\xi-x)-(il(l-2\bar k)+4k_Rk_I)(\eta-y)}q(e_{1_{\Xi t}}\tilde\mu_2^+)\nonumber\\
&-\frac{1}{2\pi}\int_{2k_R}^{0}dl\int_{-\infty}^{\infty}d\xi\int_{t}^{T}d\tau\,e^{-i(l-2k_R)(\xi-x)+(il(l-2\bar k)+4k_Rk_I)y+(\omega(-\bar k,l)+\omega(-\bar k,-2k_R))(\tau-t)}H(\xi,\tau,-\bar k,l) (e_{1_{\Xi \tau}}|_{\eta=0}\tilde \phi_2^+)\nonumber\\
&+\frac{1}{2\pi}\int_{0}^{\infty}dl\int_{-\infty}^{\infty}d\xi\int_{0}^{t}d\tau\,e^{-i(l-2k_R)(\xi-x)+(il(l-2\bar k)+4k_Rk_I)y+(\omega(-\bar k,l)+\omega(-\bar k,-2k_R))(\tau-t)}H(\xi,\tau,-\bar k,l) (e_{1_{\Xi \tau}}|_{\eta=0}\tilde \phi_2^+)\nonumber\\
&+\frac{i}{2\pi}\int_{-\infty}^{2k_R}dl\int_{-\infty}^{\infty}d\xi\int_y^\infty d\eta\, e^{-i(l-2k_R)(\xi-x)-(il(l-2\bar k)+4k_Rk_I)(\eta-y)}q(e_{1_{\Xi t}}\tilde\mu_2^+).
\end{align}}
Replacing $l$ by $l+2k_R$, we find:
{\small\begin{align}\label{tildemuplus2b}
&e_{1_{Xt}}\tilde \mu_2^+=e_{1_{Xt}}-\frac{i}{2\pi}\int_{0}^{\infty}dl\int_{-\infty}^{\infty}d\xi\int_0^y d\eta\, e^{-il(\xi-x)-il(l+2 k)(\eta-y)}q(e_{1_{\Xi t}}\tilde\mu_2^+)\nonumber\\
&-\frac{1}{2\pi}\int_0^{-2k_R}dl\int_{-\infty}^{\infty}d\xi\int_{t}^{T}d\tau\,e^{-il(\xi-x)+il(l+2k)y+\omega(k,l)(\tau-t)}H(\xi,\tau, k,l) (e_{1_{\Xi \tau}}|_{\eta=0}\tilde \phi_2^+)\nonumber\\
&+\frac{1}{2\pi}\int_{-2k_R}^{\infty}dl\int_{-\infty}^{\infty}d\xi\int_{0}^{t}d\tau\,e^{-il(\xi-x)+il(l+2k)y+\omega(k,l)(\tau-t)}H(\xi,\tau, k,l) (e_{1_{\Xi \tau}}|_{\eta=0}\tilde \phi_2^+)\nonumber\\
&+\frac{i}{2\pi}\int_{-\infty}^{0}dl\int_{-\infty}^{\infty}d\xi\int_y^\infty d\eta\, e^{-il(\xi-x)-il(l+2 k)(\eta-y)}q(e_{1_{\Xi t}}\tilde\mu_2^+),\quad k_R\leq 0,\, k_I\geq 0.
\end{align}}
Multiplying this equation by $f_1^+$, we obtain an equation which the same as the equation \eqref{dplus1i} satisfied by the d-bar derivative of $\mu_1^+$. Hence, by uniqueness we find equation \eqref{dbarplus1i}. A similar derivation yields equation \eqref{dbarminus2i}.

In the case of the d-bar derivative of $\mu_2^+$, equation \eqref{plus2i} implies that
\begin{align}\label{dplus2i}
\frac{\partial\mu_2^+}{\partial \bar k}&=-\frac{i}{2\pi}\int_{0}^{\infty}dl\int_{-\infty}^{\infty}d\xi\int_{0}^{y}d\eta\, e^{-il(\xi-x)-il(l+2k)(\eta-y)}q\frac{\partial\mu_2^+}{\partial \bar k}\nonumber\\
&+\frac{1}{2\pi}\int_{0}^{\infty}dl\int_{-\infty}^{\infty}d\xi\int_{0}^{t}d\tau\,e^{-il(\xi-x)+il(l+2k)y+\omega(k,l)(\tau-t)}H\frac{\partial\phi_2^+}{\partial \bar k}\nonumber\\
&+\frac{i}{2\pi}\int_{-\infty}^{0}dl\int_{-\infty}^{\infty}d\xi\int_{y}^{\infty}d\eta\,e^{-il(\xi-x)-il(l+2k)(\eta-y)}q\frac{\partial\mu_2^+}{\partial \bar k}, \quad k_R\geq 0,\, k_I\geq 0.
\end{align}
Thus,
\begin{align}\label{dplus2i2}
\frac{\partial\mu_2^+}{\partial \bar k}&=-\frac{i}{2\pi}\int_{0}^{\infty}dl\int_{-\infty}^{\infty}d\xi\int_{0}^{y}d\eta\, e^{-il(\xi-x)-il(l+2k)(\eta-y)}q\frac{\partial\mu_2^+}{\partial \bar k}\nonumber\\
&+\frac{1}{(2\pi)^2}\int_{0}^{\infty}dl\int_{-\infty}^{\infty}d\xi\int_{0}^{\infty}d\eta \,\delta(\eta)\int_{0}^{t}d\tau\,e^{-il(\xi-x)+il(l+2k)y+\omega(k,l)(\tau-t)}H\frac{\partial\mu_2^+}{\partial \bar k}\nonumber\\
&+\frac{i}{2\pi}\int_{-\infty}^{0}dl\int_{-\infty}^{\infty}d\xi\int_{y}^{\infty}d\eta\,e^{-il(\xi-x)-il(l+2k)(\eta-y)}q\frac{\partial\mu_2^+}{\partial \bar k}, \quad k_R\geq 0,\, k_I\geq 0.
\end{align}
This is a \textit{homogeneous} equation for $\partial \mu_2^+/\partial \bar k$, hence it admits the trivial zero solution. Therefore, uniqueness implies equation \eqref{dbarplus2i}. Equation \eqref{dbarminus1i} can be obtained in an analogous way.
\QED
\end{PROOF}

Now, we will compute the jumps of $\mu$ across the real and the imaginary axis. Subtracting \eqref{plus2i} from \eqref{plus1i} and \eqref{minus2i} from \eqref{minus1i}, evaluating the resulting expressions at $k_R=0$ and using the notation introduced in \eqref{deltamu}, we find the following equations for $\delta\mu^\pm$:
\begin{align}\label{deltaplus}
\delta\mu^+&=-\frac{i}{2\pi}\int_{0}^{\infty}dl\int_{-\infty}^{\infty}d\xi\int_{0}^{y}d\eta\, e^{-il(\xi-x)-il(l+2ik_I)(\eta-y)}q\,\delta\mu^+\nonumber\\
&+\frac{1}{2\pi}\int_{0}^{\infty}dl\int_{-\infty}^{\infty}d\xi\int_{0}^{t}d\tau\,e^{-il(\xi-x)+il(l+2ik_I)y+\omega(ik_I,l)(\tau-t)}H(\xi,\tau,ik_I,l)\,\delta\phi^+\nonumber\\
&+\frac{i}{2\pi}\int_{-\infty}^{0}dl\int_{-\infty}^{\infty}d\xi\int_{y}^{\infty}d\eta\,e^{-il(\xi-x)-il(l+2ik_I)(\eta-y)}q\,\delta\mu^+, \quad k_I\geq 0,
\end{align}
and
\begin{align}\label{deltaminus}
\delta\mu^-&=\frac{i}{2\pi}\int_{0}^{\infty}dl\int_{-\infty}^{\infty}d\xi\int_{y}^{\infty}d\eta\, e^{-il(\xi-x)-il(l+2ik_I)(\eta-y)}q\,\delta\mu^-\nonumber\\
&-\frac{1}{2\pi}\int_{-\infty}^0dl\int_{-\infty}^{\infty}d\xi\int_{t}^{T}d\tau\,e^{-il(\xi-x)+il(l+2ik_I)y+\omega(ik_I,l)(\tau-t)}H(\xi,\tau,ik_I,l)\,\delta\phi^-\nonumber\\
&-\frac{i}{2\pi}\int_{-\infty}^{0}dl\int_{-\infty}^{\infty}d\xi\int_{0}^{y}d\eta\,e^{-il(\xi-x)-il(l+2ik_I)(\eta-y)}q\,\delta\mu^-, \quad k_I\leq 0.
\end{align}
Evaluating \eqref{deltaplus} at $t=0$, we have:
\begin{align}\label{delta0plus}
\delta\mu^+_0&=-\frac{i}{2\pi}\int_{0}^{\infty}dl\int_{-\infty}^{\infty}d\xi\int_{0}^{y}d\eta\, e^{-il(\xi-x)-il(l+2ik_I)(\eta-y)}q_0\,\delta\mu^+_0\nonumber\\
&+\frac{i}{2\pi}\int_{-\infty}^{0}dl\int_{-\infty}^{\infty}d\xi\int_{y}^{\infty}d\eta\,e^{-il(\xi-x)-il(l+2ik_I)(\eta-y)}q_0\,\delta\mu^+_0, 
\end{align}
hence by uniqueness $\delta\mu_0^+=0$. Evaluating \eqref{deltaplus} at $y=0$, we find:
\begin{align}\label{deltaplus0}
\delta\phi^+&=\frac{1}{2\pi}\int_{0}^{\infty}dl\int_{-\infty}^{\infty}d\xi\int_{0}^{t}d\tau\,e^{-il(\xi-x)+\omega(ik_I,l)(\tau-t)}H(\xi,\tau,ik_I,l)\,\delta\phi^+\nonumber\\
&+\frac{i}{2\pi}\int_{-\infty}^{0}dl\int_{-\infty}^{\infty}d\xi\int_{0}^{\infty}d\eta\,e^{-il(\xi-x)-il(l+2ik_I)\eta}q\,\delta\mu^+.
\end{align}
Then, using the global relation \eqref{gri}, we find:
\begin{align}\label{deltaplus02}
\delta\phi^+&=\frac{1}{2\pi}\int_{0}^{\infty}dl\int_{-\infty}^{\infty}d\xi\int_{0}^{t}d\tau\,e^{-il(\xi-x)+\omega(ik_I,l)(\tau-t)}H(\xi,\tau,ik_I,l)\,\delta\phi^+\nonumber\\
&+\frac{1}{2\pi}\int_{-\infty}^0dl\int_{-\infty}^{\infty}d\xi\int_{0}^{t}d\tau\,e^{-il(\xi-x)+\omega(ik_I,l)(\tau-t)}H(\xi,\tau,ik_I,l)\,\delta\phi^+\nonumber\\
&-\frac{i}{2\pi}\int_{0}^{\infty}dl\int_{-\infty}^{\infty}d\xi\int_{0}^{\infty}d\eta\,e^{-il(\xi-x)-il(l+2ik_I)\eta-\omega(ik_I,l)t}q_0\,\delta\mu_0^+.
\end{align}
Since $\delta\mu_0^+=0$, uniqueness implies that $\delta \phi^+=0$. Therefore, the equation \eqref{deltaplus} for $\delta \mu^+$ becomes homogeneous, hence
\begin{equation}\label{plusdelta}
\delta \mu^+(x,y,t,k_I)=0.
\end{equation}
Similarly, we can show that
\begin{equation}\label{minusdelta}
\delta \mu^-(x,y,t,k_I)=0.
\end{equation}
Define the functions $\Delta \mu_1$ and $\Delta \mu_2$ by equations \eqref{Deltamu}. The following remark is of crucial importance to the subsequent derivations.
\vskip 2mm
\noindent\textbf{Remark 4.1.} For $k_I=0$, the exponentials appearing in equations \eqref{plus1i}-\eqref{plus2i} possess \textit{purely imaginary} exponents, hence the relevant integrals remain bounded \textit{regardless of the choice of the limits of integration}. 
\vskip 2mm
\noindent Evaluating equations \eqref{plus1i} and \eqref{minus1i} at $k_I=0$, we have:
\begin{align}\label{Dplus1i}
&\mu_1^+\Big |_{k_I=0}=1-\frac{i}{2\pi}\int_{0}^{\infty}dl\int_{-\infty}^{\infty}d\xi\int_{0}^{y}d\eta\, e^{-il(\xi-x)-il(l+2k_R)(\eta-y)}q\mu_1^+\Big |_{k_I=0}\nonumber\\
&-\frac{1}{2\pi}\int_{0}^{-2k_R}dl\int_{-\infty}^{\infty}d\xi\int_{t}^Td\tau\,e^{-il(\xi-x)+il(l+2k_R)y+\omega(k_R,l)(\tau-t)}H(\xi, \tau,k_R,l)\phi_1^+\Big |_{k_I=0}\nonumber\\
&+\frac{1}{2\pi}\int_{-2k_R}^{\infty}dl\int_{-\infty}^{\infty}d\xi\int_{0}^{t}d\tau\,e^{-il(\xi-x)+il(l+2k_R)y+\omega(k_R,l)(\tau-t)}H(\xi,\tau,k_R,l)\phi_1^+\Big |_{k_I=0}\nonumber\\
&+\frac{i}{2\pi}\int_{-\infty}^{0}dl\int_{-\infty}^{\infty}d\xi\int_{y}^{\infty}d\eta\,e^{-il(\xi-x)-il(l+2k_R)(\eta-y)}q\mu_1^+\Big |_{k_I=0}, \quad k_R\leq 0
\end{align}
and
\begin{align}\label{Dminus1initiali}
&\mu_1^-\Big |_{k_I=0}=1+\frac{i}{2\pi}\int_{0}^{\infty}dl\int_{-\infty}^{\infty}d\xi\int_{y}^\infty d\eta\, e^{-il(\xi-x)-il(l+2k_R)(\eta-y)}q\mu_1^-\Big |_{k_I=0}\nonumber\\
&-\frac{1}{2\pi}\int_{-\infty}^0dl\int_{-\infty}^{\infty}d\xi\int_{t}^{T}d\tau\,e^{-il(\xi-x)+il(l+2k_R)y+\omega(k_R,l)(\tau-t)}H(\xi,\tau,k_R,l)\phi_1^-\Big |_{k_I=0}\nonumber\\
&-\frac{i}{2\pi}\int_{-\infty}^{0}dl\int_{-\infty}^{\infty}d\xi\int_{0}^{y}d\eta\,e^{-il(\xi-x)-il(l+2k)(\eta-y)}q\mu_1^-\Big |_{k_I=0},\quad  k_R\leq 0.
\end{align}
Adding and subtracting the following term to equations \eqref{Dplus1i} and \eqref{Dminus1initiali},
\begin{align*}
&-\frac{1}{2\pi}\int_{0}^{-2k_R}dl\int_{-\infty}^{\infty}d\xi\int_{t}^Td\tau\,e^{-il(\xi-x)+il(l+2k_R)y+\omega(k_R,l)(\tau-t)}H(\xi, \tau,k_R,l)\phi_1^-\Big |_{k_I=0}\nonumber\\
&+\frac{1}{2\pi}\int_{-2k_R}^{\infty}dl\int_{-\infty}^{\infty}d\xi\int_{0}^{t}d\tau\,e^{-il(\xi-x)+il(l+2k_R)y+\omega(k_R,l)(\tau-t)}H(\xi,\tau,k_R,l)\phi_1^-\Big |_{k_I=0}, 
\end{align*}
they yield the following expressions:
\begin{align}\label{Dminus1i}
&\mu_1^-\Big |_{k_I=0}=1+\frac{i}{2\pi}\int_{0}^{\infty}dl\int_{-\infty}^{\infty}d\xi\int_{y}^\infty d\eta\, e^{-il(\xi-x)-il(l+2k_R)(\eta-y)}q\mu_1^-\Big |_{k_I=0}\nonumber\\
&-\frac{1}{2\pi}\int_{0}^{-2k_R}dl\int_{-\infty}^{\infty}d\xi\int_{t}^Td\tau\,e^{-il(\xi-x)+il(l+2k_R)y+\omega(k_R,l)(\tau-t)}H(\xi, \tau,k_R,l)\phi_1^-\Big |_{k_I=0}\nonumber\\
&+\frac{1}{2\pi}\int_{-2k_R}^{\infty}dl\int_{-\infty}^{\infty}d\xi\int_{0}^{t}d\tau\,e^{-il(\xi-x)+il(l+2k_R)y+\omega(k_R,l)(\tau-t)}H(\xi,\tau,k_R,l)\phi_1^-\Big |_{k_I=0}\nonumber\\
&-\frac{i}{2\pi}\int_{-\infty}^{0}dl\int_{-\infty}^{\infty}d\xi\int_{0}^{y}d\eta\,e^{-il(\xi-x)-il(l+2k_R)(\eta-y)}q\mu_1^-\Big |_{k_I=0}\nonumber\\
&-\frac{1}{2\pi}\int_{-\infty}^0dl\int_{-\infty}^{\infty}d\xi\int_{t}^{T}d\tau\,e^{-il(\xi-x)+il(l+2k_R)y+\omega(k_R,l)(\tau-t)}H(\xi,\tau,k_R,l)\phi_1^-\Big |_{k_I=0}\nonumber\\
&+\frac{1}{2\pi}\int_{0}^{-2k_R}dl\int_{-\infty}^{\infty}d\xi\int_{t}^{T}d\tau\,e^{-il(\xi-x)+il(l+2k_R)y+\omega(k_R,l)(\tau-t)}H(\xi,\tau,k_R,l)\phi_1^-\Big |_{k_I=0}\nonumber\\
&-\frac{1}{2\pi}\int_{-2k_R}^{\infty}dl\int_{-\infty}^{\infty}d\xi\int_{0}^{t}d\tau\,e^{-il(\xi-x)+il(l+2k_R)y+\omega(k_R,l)(\tau-t)}H(\xi,\tau,k_R,l)\phi_1^-\Big |_{k_I=0},\quad  k_R\leq 0.
\end{align}
By rearranging the integrals with respect to $\tau$ appropriately, we can write the last three terms on the RHS of \eqref{Dminus1i} as
\begin{align}\label{minus1iterms}
&-\frac{1}{2\pi}\int_{-\infty}^0dl\int_{-\infty}^{\infty}d\xi\int_{0}^{T}d\tau\,e^{-il(\xi-x)+il(l+2k_R)y+\omega(k_R,l)(\tau-t)}H(\xi,\tau,k_R,l)\phi_1^-\Big |_{k_I=0}\nonumber\\
&+\frac{1}{2\pi}\int_{-\infty}^0dl\int_{-\infty}^{\infty}d\xi\int_{0}^{t}d\tau\,e^{-il(\xi-x)+il(l+2k_R)y+\omega(k_R,l)(\tau-t)}H(\xi,\tau,k_R,l)\phi_1^-\Big |_{k_I=0}\nonumber\\
&+\frac{1}{2\pi}\int_{0}^{-2k_R}dl\int_{-\infty}^{\infty}d\xi\int_{0}^{T}d\tau\,e^{-il(\xi-x)+il(l+2k_R)y+\omega(k_R,l)(\tau-t)}H(\xi,\tau,k_R,l)\phi_1^-\Big |_{k_I=0}\nonumber\\
&-\frac{1}{2\pi}\int_{0}^{\infty}dl\int_{-\infty}^{\infty}d\xi\int_{0}^{t}d\tau\,e^{-il(\xi-x)+il(l+2k_R)y+\omega(k_R,l)(\tau-t)}H(\xi,\tau,k_R,l)\phi_1^-\Big |_{k_I=0}.
\end{align}
Furthermore, by employing the global relation \eqref{gri}, these terms become:
\begin{align}\label{minus1iterms2}
&-\frac{1}{2\pi}\int_{-\infty}^0dl\int_{-\infty}^{\infty}d\xi\int_{0}^{T}d\tau\,e^{-il(\xi-x)+il(l+2k_R)y+\omega(k_R,l)(\tau-t)}H(\xi,\tau,k_R,l)\phi_1^-\Big |_{k_I=0}\nonumber\\
&+\frac{1}{2\pi}\int_{0}^{-2k_R}dl\int_{-\infty}^{\infty}d\xi\int_{0}^{T}d\tau\,e^{-il(\xi-x)+il(l+2k_R)y+\omega(k_R,l)(\tau-t)}H(\xi,\tau,k_R,l)\phi_1^-\Big |_{k_I=0}\nonumber\\
&+\frac{i}{2\pi}\int_{-\infty}^0 dl\int_{-\infty}^{\infty}d\xi\int_{0}^{\infty}d\eta\,e^{-il(\xi-x)-il(l+2k_R)(\eta-y)}q\mu_1^-\Big |_{k_I=0}\nonumber\\
&-\frac{i}{2\pi}\int_{-\infty}^0 dl\int_{-\infty}^{\infty}d\xi\int_{0}^{\infty}d\eta\,e^{-il(\xi-x)-il(l+2k_R)(\eta-y)-\omega(k_R,l)t}q_0\mu_1^-\Big |_{k_I=0,\,t=0}\nonumber\\
&-\frac{i}{2\pi}\int_{0}^\infty dl\int_{-\infty}^{\infty}d\xi\int_{0}^{\infty}d\eta\,e^{-il(\xi-x)-il(l+2k_R)(\eta-y)}q\mu_1^-\Big |_{k_I=0}\nonumber\\
&+\frac{i}{2\pi}\int_0^{\infty} dl\int_{-\infty}^{\infty}d\xi\int_{0}^{\infty}d\eta\,e^{-il(\xi-x)-il(l+2k_R)(\eta-y)-\omega(k_R,l)t}q_0\mu_1^-\Big |_{k_I=0,\,t=0}.
\end{align}
Hence, we find the following equation for $\Delta \mu_1$:
\begin{align}\label{Delta1i}
&\Delta \mu_1=-\frac{i}{2\pi}\int_{0}^{\infty}dl\int_{-\infty}^{\infty}d\xi\int_0^{y} d\eta\, e^{-il(\xi-x)-il(l+2k_R)(\eta-y)}q\Delta \mu_1\nonumber\\
&-\frac{1}{2\pi}\int_{0}^{-2k_R}dl\int_{-\infty}^{\infty}d\xi\int_{t}^Td\tau\,e^{-il(\xi-x)+il(l+2k_R)y+\omega(k_R,l)(\tau-t)}H(\xi, \tau,k_R,l)\Delta \phi_1\nonumber\\
&+\frac{1}{2\pi}\int_{-2k_R}^{\infty}dl\int_{-\infty}^{\infty}d\xi\int_{0}^{t}d\tau\,e^{-il(\xi-x)+il(l+2k_R)y+\omega(k_R,l)(\tau-t)}H(\xi,\tau,k_R,l)\Delta \phi_1\nonumber\\
&+\frac{i}{2\pi}\int_{-\infty}^{0}dl\int_{-\infty}^{\infty}d\xi\int_{y}^\infty d\eta\,e^{-il(\xi-x)-il(l+2k_R)(\eta-y)}q\Delta \mu_1\nonumber\\
&-\frac{i}{2\pi}\int_0^{\infty} dl\int_{-\infty}^{\infty}d\xi\int_{0}^{\infty}d\eta\,e^{-il(\xi-x)-il(l+2k_R)(\eta-y)-\omega(k_R,l)t}q_0\mu_1^-\Big |_{k_I=0,\,t=0}\nonumber\\
&-\frac{1}{2\pi}\int_{0}^{-2k_R}dl\int_{-\infty}^{\infty}d\xi\int_{0}^{T}d\tau\,e^{-il(\xi-x)+il(l+2k_R)y+\omega(k_R,l)(\tau-t)}H(\xi,\tau,k_R,l)\phi_1^-\Big|_{k_I=0}\nonumber\\
&+\frac{1}{2\pi}\int_{-\infty}^0dl\int_{-\infty}^{\infty}d\xi\int_{0}^{T}d\tau\,e^{-il(\xi-x)+il(l+2k_R)y+\omega(k_R,l)(\tau-t)}H(\xi,\tau,k_R,l)\phi_1^-\Big |_{k_I=0}\nonumber\\
&+\frac{i}{2\pi}\int_{-\infty}^0 dl\int_{-\infty}^{\infty}d\xi\int_{0}^{\infty}d\eta\,e^{-il(\xi-x)-il(l+2k_R)(\eta-y)-\omega(k_R,l)t}q_0\mu_1^-\Big |_{k_I=0,\,t=0},\quad  k_R\leq 0.
\end{align}
Let us now introduce the notation
\begin{equation}\label{eXt2}
e_{2_{Xt}}:=e_{Xt}(x,y,t,k_R,\lambda)=e^{-i(\lambda+2k_R)x+i\lambda(\lambda+2k_R)y+\omega(k_R,\lambda)t-8ik_R^3 t}.
\end{equation}
Then, under the change of variables 
\begin{equation}\label{transD1}
k_R\mapsto -k_R-\frac{\lambda}{2}, \quad k_I\mapsto i \frac{\lambda}{2}, \quad l\mapsto l+2k_R+\lambda, \quad \lambda \in \mathbb R
\end{equation}
and introducing the notation 
\begin{equation}\label{muplus1breve}
 \breve \mu_j^+(x,y,t,k_R,k_I)=\mu_j^+(x,y,t,-k_R-\frac{\lambda}{2},i \frac{\lambda}{2}), \quad j=1,2,
\end{equation}
equation \eqref{plus1i} takes the following form:
\begin{align}\label{plus1itrans}
&e_{2_{Xt}}\breve\mu_1^+=e_{2_{Xt}}-\frac{i}{2\pi}\int_{-2k_R-\lambda}^{\infty}dl\int_{-\infty}^{\infty}d\xi\int_{0}^{y}d\eta\, e^{-il(\xi-x)-il(l+2k_R)(\eta-y)}q\left(e_{2_{\Xi t}}\breve\mu_1^+\right)\nonumber\\
&-\frac{1}{2\pi}\int_{-2k_R-\lambda}^0 dl\int_{-\infty}^{\infty}d\xi\int_{t}^T d\tau\,e^{-il(\xi-x)+il(l+2k_R)y+\omega(k_R,l)(\tau-t)}H\left(e_{2_{\Xi \tau}}|_{\eta=0}\breve\phi_1^+\right)\nonumber\\
&+\frac{1}{2\pi}\int_{0}^{\infty}dl\int_{-\infty}^{\infty}d\xi\int_{0}^{t}d\tau\,e^{-il(\xi-x)+il(l+2k_R)y+\omega(k_R,l)(\tau-t)}H \left(e_{2_{\Xi \tau}}|_{\eta=0}\breve\phi_1^+\right)\nonumber\\
&+\frac{i}{2\pi}\int_{-\infty}^{-2k_R-\lambda}dl\int_{-\infty}^{\infty}d\xi\int_{y}^{\infty}d\eta\,e^{-il(\xi-x)-il(l+2k_R)(\eta-y)}q\left(e_{2_{\Xi t}}\breve\mu_1^+\right), \quad k_R\leq 0, \,\lambda\geq -2k_R.
\end{align}
\vskip 2mm
\noindent \textbf{Remark 4.2.} Note that the function $H$, defined by \eqref{H}, remains \textit{invariant} under the change of variables \eqref{transD1}.
\vskip 2mm
\noindent Rearranging the integrals with respect to $l$, we have:
\begin{align}\label{plus1itrans2}
&e_{2_{Xt}}\breve\mu_1^+=e_{2_{Xt}}-\frac{i}{2\pi}\int_{0}^{\infty}dl\int_{-\infty}^{\infty}d\xi\int_{0}^{y}d\eta\, e^{-il(\xi-x)-il(l+2k_R)(\eta-y)}q\left(e_{2_{\Xi t}}\breve\mu_1^+\right)\nonumber\\
&-\frac{1}{2\pi}\int_{0}^{-2k_R} dl\int_{-\infty}^{\infty}d\xi\int_{t}^T d\tau\,e^{-il(\xi-x)+il(l+2k_R)y+\omega(k_R,l)(\tau-t)}H\left(e_{2_{\Xi \tau}}|_{\eta=0}\breve\phi_1^+\right)\nonumber\\
&+\frac{1}{2\pi}\int_{-2k_R}^{\infty}dl\int_{-\infty}^{\infty}d\xi\int_{0}^{t}d\tau\,e^{-il(\xi-x)+il(l+2k_R)y+\omega(k_R,l)(\tau-t)}H \left(e_{2_{\Xi \tau}}|_{\eta=0}\breve\phi_1^+\right)\nonumber\\
&+\frac{i}{2\pi}\int_{-\infty}^{0}dl\int_{-\infty}^{\infty}d\xi\int_{y}^{\infty}d\eta\,e^{-il(\xi-x)-il(l+2k_R)(\eta-y)}q\left(e_{2_{\Xi t}}\breve\mu_1^+\right),\nonumber\\
&-\frac{i}{2\pi}\int_{-2k_R-\lambda}^{0}dl\int_{-\infty}^{\infty}d\xi\int_{0}^{y}d\eta\, e^{-il(\xi-x)-il(l+2k_R)(\eta-y)}q\left(e_{2_{\Xi t}}\breve\mu_1^+\right)\nonumber\\
&-\frac{1}{2\pi}\int_{-2k_R-\lambda}^0 dl\int_{-\infty}^{\infty}d\xi\int_{t}^T d\tau\,e^{-il(\xi-x)+il(l+2k_R)y+\omega(k_R,l)(\tau-t)}H\left(e_{2_{\Xi \tau}}|_{\eta=0}\breve\phi_1^+\right)\nonumber\\
&+\frac{1}{2\pi}\int_{0}^{-2k_R}dl\int_{-\infty}^{\infty}d\xi\int_{0}^{T}d\tau\,e^{-il(\xi-x)+il(l+2k_R)y+\omega(k_R,l)(\tau-t)}H \left(e_{2_{\Xi \tau}}|_{\eta=0}\breve\phi_1^+\right)\nonumber\\
&+\frac{i}{2\pi}\int_{0}^{-2k_R-\lambda}dl\int_{-\infty}^{\infty}d\xi\int_{y}^{\infty}d\eta\,e^{-il(\xi-x)-il(l+2k_R)(\eta-y)}q\left(e_{2_{\Xi t}}\breve\mu_1^+\right), \quad k_R\leq 0, \,\lambda\geq -2k_R.
\end{align}
Furthermore, employing the global relation \eqref{gri} in the last four terms on the RHS of the above expression, we find:
\begin{align}\label{plus1itrans3}
&e_{2_{Xt}}\breve\mu_1^+=e_{2_{Xt}}-\frac{i}{2\pi}\int_{0}^{\infty}dl\int_{-\infty}^{\infty}d\xi\int_{0}^{y}d\eta\, e^{-il(\xi-x)-il(l+2k_R)(\eta-y)}q\left(e_{2_{\Xi t}}\breve\mu_1^+\right)\nonumber\\
&-\frac{1}{2\pi}\int_{0}^{-2k_R} dl\int_{-\infty}^{\infty}d\xi\int_{t}^T d\tau\,e^{-il(\xi-x)+il(l+2k_R)y+\omega(k_R,l)(\tau-t)}H\left(e_{2_{\Xi \tau}}|_{\eta=0}\breve\phi_1^+\right)\nonumber\\
&+\frac{1}{2\pi}\int_{-2k_R}^{\infty}dl\int_{-\infty}^{\infty}d\xi\int_{0}^{t}d\tau\,e^{-il(\xi-x)+il(l+2k_R)y+\omega(k_R,l)(\tau-t)}H \left(e_{2_{\Xi \tau}}|_{\eta=0}\breve\phi_1^+\right)\nonumber\\
&+\frac{i}{2\pi}\int_{-\infty}^{0}dl\int_{-\infty}^{\infty}d\xi\int_{y}^{\infty}d\eta\,e^{-il(\xi-x)-il(l+2k_R)(\eta-y)}q\left(e_{2_{\Xi t}}\breve\mu_1^+\right),\nonumber\\
&+\int_0^{-2k_R-\lambda}dl\, e^{ilx+il(l+2k_R)y-\omega(k_R,l)t}\Big[p_1^+(-k_R-\lambda,l+2k_R+\lambda)+r_1^+(-k_R-\lambda,l+2k_R+\lambda)\Big]\nonumber\\
&+\int_{0}^{-2k_R}dl\,e^{ilx+il(l+2k_R)y-\omega(k_R,l)t}p_1^+(-k_R-\lambda,l+2k_R+\lambda), \quad k_R\leq 0, \,\lambda\geq -2k_R,
\end{align}
where the functions $p_1^+$ and $r_1^+$ are defined as follows:
\begin{subequations}\label{p1r1}
\begin{align}
p_1^+(k,l)&:=\frac{1}{2\pi}\int_{-\infty}^{\infty}d\xi\int_{0}^{T}d\tau\,e^{-il\xi+\omega(k,l)\tau}H(\xi,\tau,k,l) \phi_1^+,\label{p1}\\
r_1^+(k,l)&:=\frac{i}{2\pi}\int_{-\infty}^{\infty}d\xi\int_{0}^{\infty}d\eta\,e^{-il\xi-il(l+2k)\eta}q_0\mu_1^+\Big|_{t=0}.\label{r1}
\end{align}
\end{subequations}
The change of variables \eqref{transD1} applied on equation \eqref{plus2i} yields:
\begin{align}\label{plus2itrans}
e_{2_{Xt}}\breve\mu_2^+&=e_{2_{Xt}}-\frac{i}{2\pi}\int_{-2k_R-\lambda}^{\infty}dl\int_{-\infty}^{\infty}d\xi\int_{0}^{y}d\eta\, e^{-il(\xi-x)-il(l+2k_R)(\eta-y)}q\left(e_{2_{\Xi t}}\breve\mu_2^+\right)\nonumber\\
&+\frac{1}{2\pi}\int_{-2k_R-\lambda}^\infty dl\int_{-\infty}^{\infty}d\xi\int_{0}^t d\tau\,e^{-il(\xi-x)+il(l+2k_R)y+\omega(k_R,l)(\tau-t)}H\left(e_{2_{\Xi \tau}}|_{\eta=0}\breve\phi_2^+\right)\nonumber\\
&+\frac{i}{2\pi}\int_{-\infty}^{-2k_R-\lambda}dl\int_{-\infty}^{\infty}d\xi\int_{y}^{\infty}d\eta\,e^{-il(\xi-x)-il(l+2k_R)(\eta-y)}q\left(e_{2_{\Xi t}}\breve\mu_2^+\right), \quad k_R\leq 0, \,\lambda\leq -2k_R.
\end{align}
As in the case of $\breve \mu_1^+$, rearranging the integrals with respect to $l$ yields:
\begin{align}\label{plus2itrans2}
&e_{2_{Xt}}\breve\mu_2^+=e_{2_{Xt}}-\frac{i}{2\pi}\int_{0}^{\infty}dl\int_{-\infty}^{\infty}d\xi\int_{0}^{y}d\eta\, e^{-il(\xi-x)-il(l+2k_R)(\eta-y)}q\left(e_{2_{\Xi t}}\breve\mu_2^+\right)\nonumber\\
&-\frac{1}{2\pi}\int_{0}^{-2k_R} dl\int_{-\infty}^{\infty}d\xi\int_{t}^T d\tau\,e^{-il(\xi-x)+il(l+2k_R)y+\omega(k_R,l)(\tau-t)}H\left(e_{2_{\Xi \tau}}|_{\eta=0}\breve\phi_2^+\right)\nonumber\\
&+\frac{1}{2\pi}\int_{-2k_R}^{\infty}dl\int_{-\infty}^{\infty}d\xi\int_{0}^{t}d\tau\,e^{-il(\xi-x)+il(l+2k_R)y+\omega(k_R,l)(\tau-t)}H \left(e_{2_{\Xi \tau}}|_{\eta=0}\breve\phi_2^+\right)\nonumber\\
&+\frac{i}{2\pi}\int_{-\infty}^{0}dl\int_{-\infty}^{\infty}d\xi\int_{y}^{\infty}d\eta\,e^{-il(\xi-x)-il(l+2k_R)(\eta-y)}q\left(e_{2_{\Xi t}}\breve\mu_2^+\right)\nonumber\\
&+\frac{i}{2\pi}\int^{-2k_R-\lambda}_{0}dl\int_{-\infty}^{\infty}d\xi\int_{0}^{\infty}d\eta\, e^{-il(\xi-x)-il(l+2k_R)(\eta-y)}q\left(e_{2_{\Xi t}}\breve\mu_2^+\right)\nonumber\\
&+\frac{1}{2\pi}\int_{-2k_R-\lambda}^{-2k_R} dl\int_{-\infty}^{\infty}d\xi\int_{0}^T d\tau\,e^{-il(\xi-x)+il(l+2k_R)y+\omega(k_R,l)(\tau-t)}H\left(e_{2_{\Xi \tau}}|_{\eta=0}\breve\phi_2^+\right)\nonumber\\
&-\frac{1}{2\pi}\int_{-2k_R-\lambda}^{0} dl\int_{-\infty}^{\infty}d\xi\int_{t}^T d\tau\,e^{-il(\xi-x)+il(l+2k_R)y+\omega(k_R,l)(\tau-t)}H\left(e_{2_{\Xi \tau}}|_{\eta=0}\breve\phi_2^+\right), \quad k_R\leq 0, \,\lambda\leq -2k_R.
\end{align}
Using the global relation \eqref{gri} we find the following equation for $\breve \mu_2^+$:
\begin{align}\label{plus2itrans3}
&e_{2_{Xt}}\breve\mu_2^+=e_{2_{Xt}}-\frac{i}{2\pi}\int_{0}^{\infty}dl\int_{-\infty}^{\infty}d\xi\int_{0}^{y}d\eta\, e^{-il(\xi-x)-il(l+2k_R)(\eta-y)}q\left(e_{2_{\Xi t}}\breve\mu_2^+\right)\nonumber\\
&-\frac{1}{2\pi}\int_{0}^{-2k_R} dl\int_{-\infty}^{\infty}d\xi\int_{t}^T d\tau\,e^{-il(\xi-x)+il(l+2k_R)y+\omega(k_R,l)(\tau-t)}H\left(e_{2_{\Xi \tau}}|_{\eta=0}\breve\phi_2^+\right)\nonumber\\
&+\frac{1}{2\pi}\int_{-2k_R}^{\infty}dl\int_{-\infty}^{\infty}d\xi\int_{0}^{t}d\tau\,e^{-il(\xi-x)+il(l+2k_R)y+\omega(k_R,l)(\tau-t)}H \left(e_{2_{\Xi \tau}}|_{\eta=0}\breve\phi_2^+\right)\nonumber\\
&+\frac{i}{2\pi}\int_{-\infty}^{0}dl\int_{-\infty}^{\infty}d\xi\int_{y}^{\infty}d\eta\,e^{-il(\xi-x)-il(l+2k_R)(\eta-y)}q\left(e_{2_{\Xi t}}\breve\mu_2^+\right)\nonumber\\
&+\int_{0}^{-2k_R} dl\,e^{ilx+il(l+2k_R)y-\omega(k_R,l)t}p_2^+(-k_R-\lambda,l+2k_R+\lambda)\nonumber\\
&+\int^{-2k_R-\lambda}_{0} dl\,e^{ilx+il(l+2k_R)y-\omega(k_R,l)t}r_2^+(-k_R-\lambda,l+2k_R+\lambda), \quad k_R\leq 0, \,\lambda\leq -2k_R,
\end{align}
where the functions $p_2^+$ and $r_2^+$ are defined as follows:
\begin{subequations}\label{p2r2}
\begin{align}
p_2^+(k,l)&=\frac{1}{2\pi}\int_{-\infty}^{\infty}d\xi\int_{0}^{T}d\tau\,e^{-il\xi+\omega(k,l)\tau}H(\xi,\tau,k,l) \phi_2^+,\label{p2}\\
r_2^+(k,l)&=\frac{i}{2\pi}\int_{-\infty}^{\infty}d\xi\int_{0}^{\infty}d\eta\,e^{-il\xi-il(l+2k)\eta}q_0\mu_2^+\Big|_{t=0}.\label{r2}
\end{align}
\end{subequations}
The above equations imply the following result:
\begin{prop}\label{D1}
The solution $\Delta\mu_1$ of equation \eqref{Delta1i} is given by
\begin{align}\label{solD1}
\Delta\mu_1(x,y,t,k_R)&=\int_{-\infty}^{-2k_R}d\lambda\, \chi_1(k_R,\lambda) \left(e_{2_{Xt}}\breve \mu_2^+\right)+\int_{-2k_R}^{\infty}d\lambda\, \chi_2(k_R,\lambda) \left(e_{2_{Xt}}\breve \mu_1^+\right)\nonumber\\
&+\int_{0}^{-2k_R}d\lambda\, \chi_3(k_R,\lambda) \left(e_{2_{Xt}}\breve \mu_2^+\right), \quad k_R\leq 0, 
\end{align}
provided that the functions $\chi_j,\ j=1,2,3,$ satisfy the following linear Volterra integral equations: 
\begin{subequations}\label{chi}
\begin{align}
&\chi_1(k_R,\lambda)+\int_{-\infty}^{\lambda}dl\, \chi_1(k_R,l)r_2^+(-k_R-l,l-\lambda)=-r_1^-(k_R,-\lambda-2k_R), \quad k_R\leq 0,\ \lambda \leq -2k_R,\label{chi1}\\
&\chi_2(k_R,\lambda)-\int_{\lambda}^\infty dl\, \chi_2(k_R,l)\left[p_1^+(-k_R-l,l-\lambda)+r_1^+(-k_R-l,l-\lambda)\right]=\nonumber\\
&=p_1^-(k_R,-\lambda-2k_R)+r_1^-(k_R,-\lambda-2k_R), \quad k_R\leq 0,\ \lambda \geq -2k_R,\label{chi2}\\ 
\nonumber\\
&\chi_3(k_R,\lambda)+\int_{0}^{-2k_R}dl\, \chi_3(k_R,l)p_2^+(-k_R-l,l-\lambda)\nonumber\\
&+\int_{0}^{\lambda}dl\, \chi_3(k_R,l)r_2^+(-k_R-l,l-\lambda)+\int_{-\infty}^{-2k_R}dl\, \chi_1(k_R,l)p_2^+(-k_R-l,l-\lambda)\nonumber\\
&+\int_{-2k_R}^{\infty}dl\, \chi_2(k_R,l)p_1^+(-k_R-l,l-\lambda)=-p_1^-(k_R,-\lambda-2k_R), \quad k_R\leq 0,\ 0\leq\lambda \leq -2k_R\label{chi3}, 
\end{align}
\end{subequations}
where the functions $p^\pm_j$ and $r^\pm_j$, $j=1,2,$ are defined by:
\begin{subequations}\label{pjrj}
\begin{align}
p_j^\pm(k,l)&:=\frac{1}{2\pi}\int_{-\infty}^{\infty}d\xi\int_{0}^{T}d\tau\,e^{-il\xi+\omega(k,l)\tau}H(\xi,\tau,k,l) \phi_j^\pm,\label{pj}\\
r_j^\pm(k,l)&:=\frac{i}{2\pi}\int_{-\infty}^{\infty}d\xi\int_{0}^{\infty}d\eta\,e^{-il\xi-il(l+2k)\eta}q_0\mu_j^\pm\Big|_{t=0}.\label{rj}
\end{align}
\end{subequations}
\end{prop}

\begin{PROOF}
Employing equations \eqref{plus1itrans3} and \eqref{plus2itrans3}, we notice that the equation \eqref{solD1} for $\mathcal S_1$ is precisely the equation \eqref{Delta1i} for $\Delta \mu_1$ as long as the functions $\chi_j,\ j=1,2,3,$ are chosen so that the forcing terms of the two equations match. The forcing term of \eqref{solD1} is equal to
\begin{align}\label{fS1}
&\int_{-\infty}^{-2k_R}d\lambda\, \chi_1(k_R,\lambda) \left[e_{2_{Xt}}+\int_{0}^{-2k_R} dl\,e^{ilx+il(l+2k_R)y-\omega(k_R,l)t}p_2^+(-k_R-\lambda,l+2k_R+\lambda)\right.\nonumber\\
&\left.+\int^{-2k_R-\lambda}_{0} dl\,e^{ilx+il(l+2k_R)y-\omega(k_R,l)t}r_2^+(-k_R-\lambda,l+2k_R+\lambda)\right]\nonumber\\
&+\int_{-2k_R}^\infty d\lambda\, \chi_2(k_R,\lambda)\left[e_{2_{Xt}}+\int_{0}^{-2k_R}dl\,e^{ilx+il(l+2k_R)y-\omega(k_R,l)t}p_1^+(-k_R-\lambda,l+2k_R+\lambda)\right.\nonumber\\
&\left.+\int_0^{-2k_R-\lambda}dl\, e^{ilx+il(l+2k_R)y-\omega(k_R,l)t}\Big[p_1^+(-k_R-\lambda,l+2k_R+\lambda)+r_1^+(-k_R-\lambda,l+2k_R+\lambda)\Big]\right]\nonumber\\
&+\int_{0}^{-2k_R}d\lambda\, \chi_3(k_R,\lambda) \left[e_{2_{Xt}}+\int_{0}^{-2k_R} dl\,e^{ilx+il(l+2k_R)y-\omega(k_R,l)t}p_2^+(-k_R-\lambda,l+2k_R+\lambda)\right.\nonumber\\
&\left.+\int^{-2k_R-\lambda}_{0} dl\,e^{ilx+il(l+2k_R)y-\omega(k_R,l)t}r_2^+(-k_R-\lambda,l+2k_R+\lambda)\right].
\end{align}
In fact, interchanging the order of integration between $\lambda$ and $l$ the forcing becomes
\begin{align}\label{fS12}
&\int_{-\infty}^{-2k_R}d\lambda\, \chi_1(k_R,\lambda)e_{2_{Xt}}+\int_{0}^{-2k_R}dl\, E_{Xt}\int_{-\infty}^{-2k_R} d\lambda\, \chi_1(k_R,\lambda)\,p_2^+(-k_R-\lambda,l+2k_R+\lambda)\nonumber\\
&+\int_{0}^{\infty} dl\,E_{Xt} \int_{-\infty}^{-2k_R-l}d\lambda\,\chi_1(k_R,\lambda)\,r_2^+(-k_R-\lambda,l+2k_R+\lambda)\nonumber\\
&+\int_{-2k_R}^{\infty}d\lambda\, \chi_2(k_R,\lambda)e_{2_{Xt}}+\int_{0}^{-2k_R}dl\, E_{Xt}\int_{-2k_R}^\infty d\lambda\, \chi_2(k_R,\lambda)\,p_1^+(-k_R-\lambda,l+2k_R+\lambda)\nonumber\\
&-\int_{-\infty}^{0} dl\,E_{Xt} \int_{-2k_R-l}^\infty d\lambda\,\chi_2(k_R,\lambda)\left[p_1^+(-k_R-\lambda,l+2k_R+\lambda)+r_1^+(-k_R-\lambda,l+2k_R+\lambda)\right]\nonumber\\
&+\int_{0}^{-2k_R}d\lambda\, \chi_3(k_R,\lambda)e_{2_{Xt}}+\int_{0}^{-2k_R}dl\, E_{Xt}\int_{0}^{-2k_R} d\lambda\, \chi_3(k_R,\lambda)\,p_2^+(-k_R-\lambda,l+2k_R+\lambda)\nonumber\\
&+\int_{0}^{-2k_R} dl\,E_{Xt} \int_{0}^{-2k_R-l}d\lambda\,\chi_3(k_R,\lambda)\,r_2^+(-k_R-\lambda,l+2k_R+\lambda),
\end{align}
where the notation $E_{Xt}$ stands for the exponential
\begin{equation}\label{EXt}
E_{Xt}(x,y,t,k_R,l):=e^{i l x+i l (l+2k_R)y-\omega(k_R,l)t}.
\end{equation}
The definition \eqref{eXt2} of $e_{2_{Xt}}$ implies:
\begin{equation}\label{Ee}
E_{Xt}(x,y,t,k_R,l)=e_{2_{Xt}}(x,y,t,k_R,-2k_R-l),
\end{equation}
thus, letting $l\mapsto -2k_R-l$ in \eqref{fS12} yields the final expression for the forcing of \eqref{solD1} as
\begin{align}\label{fS13}
&\int_{-\infty}^{-2k_R}dl\, \chi_1(k_R,l)e_{2_{Xt}}+\int_{0}^{-2k_R}dl\, e_{2_{Xt}}\int_{-\infty}^{-2k_R} d\lambda\, \chi_1(k_R,\lambda)\,p_2^+(-k_R-\lambda,\lambda-l)\nonumber\\
&+\int_{-\infty}^{-2k_R} dl\,e_{2_{Xt}} \int_{-\infty}^{l}d\lambda\,\chi_1(k_R,\lambda)\,r_2^+(-k_R-\lambda,\lambda-l)\nonumber\\
&+\int_{-2k_R}^{\infty}dl\, \chi_2(k_R,l)e_{2_{Xt}}+\int_{0}^{-2k_R}dl\, e_{2_{Xt}}\int_{-2k_R}^\infty d\lambda\, \chi_2(k_R,\lambda)\,p_1^+(-k_R-\lambda,\lambda-l)\nonumber\\
&-\int_{-2k_R}^{\infty} dl\,e_{2_{Xt}} \int_{l}^\infty d\lambda\,\chi_2(k_R,\lambda)\left[p_1^+(-k_R-\lambda,\lambda-l)+r_1^+(-k_R-\lambda,\lambda-l)\right]\nonumber\\
&+\int_{0}^{-2k_R}dl\, \chi_3(k_R,l)e_{2_{Xt}}+\int_{0}^{-2k_R}dl\, e_{2_{Xt}}\int_{0}^{-2k_R} d\lambda\, \chi_3(k_R,\lambda)\,p_2^+(-k_R-\lambda,\lambda-l)\nonumber\\
&+\int_{0}^{-2k_R} dl\,e_{2_{Xt}} \int_{0}^{l}d\lambda\,\chi_3(k_R,\lambda)\,r_2^+(-k_R-\lambda,\lambda-l).
\end{align}
On the other hand, under the transformation $l\mapsto -2k_R-l$ the forcing of equation \eqref{Delta1i} becomes
\begin{align}\label{fD1}
&\int_{-2k_R}^{\infty}dl\, e_{2_{Xt}}\, p_1^-(k_R,-2k_R-l)-\int_{-\infty}^{-2k_R}dl\, e_{2_{Xt}}\, r_1^-(k_R,-2k_R-l)\nonumber\\
&-\int_0^{-2k_R}dl\, e_{2_{Xt}}\, p_1^-(k_R,-2k_R-l)+\int_{-2k_R}^{\infty}dl\, e_{2_{Xt}}\, r_1^-(k_R,-2k_R-l),
\end{align}
where $p_1^-$ and $r_1^-$ are defined by equations \eqref{pj} and \eqref{rj} respectively.

Therefore, equating the two forcing terms, namely \eqref{fS13} and \eqref{fD1}, we obtain the equations \eqref{chi} for $\chi_j, \ j=1,2,3$.
\QED
\end{PROOF}

Regarding the analogous expression for $\Delta \mu_2$, we evaluate equations \eqref{minus2i} and \eqref{plus2i} at $k_I=0$ and then subtract them to obtain the equation:
\begin{align}\label{Delta2i}
\Delta \mu_2&=-\frac{i}{2\pi}\int_{0}^{\infty}dl\int_{-\infty}^{\infty}d\xi\int_0^{y} d\eta\, e^{-il(\xi-x)-il(l+2k_R)(\eta-y)}q\Delta \mu_2\nonumber\\
&+\frac{1}{2\pi}\int_{0}^{\infty}dl\int_{-\infty}^{\infty}d\xi\int_{0}^{t}d\tau\,e^{-il(\xi-x)+il(l+2k_R)y+\omega(k_R,l)(\tau-t)}H(\xi,\tau,k_R,l)\Delta \phi_2\nonumber\\
&+\frac{i}{2\pi}\int_{-\infty}^{0}dl\int_{-\infty}^{\infty}d\xi\int_{y}^\infty d\eta\,e^{-il(\xi-x)-il(l+2k_R)(\eta-y)}q\Delta \mu_2\nonumber\\
&-\frac{i}{2\pi}\int_0^{\infty} dl\int_{-\infty}^{\infty}d\xi\int_{0}^{\infty}d\eta\,e^{-il(\xi-x)-il(l+2k_R)(\eta-y)-\omega(k_R,l)t}q\mu_2^-\Big |_{k_I=0}\nonumber\\
&+\frac{1}{2\pi}\int_{0}^{\infty}dl\int_{-\infty}^{\infty}d\xi\int_{0}^{t}d\tau\,e^{-il(\xi-x)+il(l+2k_R)y+\omega(k_R,l)(\tau-t)}H(\xi,\tau,k_R,l)\phi_2^-\Big|_{k_I=0}\nonumber\\
&+\frac{1}{2\pi}\int_{-\infty}^{-2k_R}dl\int_{-\infty}^{\infty}d\xi\int_{t}^{T}d\tau\,e^{-il(\xi-x)+il(l+2k_R)y+\omega(k_R,l)(\tau-t)}H(\xi,\tau,k_R,l)\phi_2^-\Big |_{k_I=0}\nonumber\\
&-\frac{1}{2\pi}\int_{-2k_R}^0dl\int_{-\infty}^{\infty}d\xi\int_{0}^{t}d\tau\,e^{-il(\xi-x)+il(l+2k_R)y+\omega(k_R,l)(\tau-t)}H(\xi,\tau,k_R,l)\phi_2^-\Big |_{k_I=0}\nonumber\\
&+\frac{i}{2\pi}\int_{-\infty}^0 dl\int_{-\infty}^{\infty}d\xi\int_{0}^{\infty}d\eta\,e^{-il(\xi-x)-il(l+2k_R)(\eta-y)-\omega(k_R,l)t}q\mu_2^-\Big |_{k_I=0},\quad  k_R\geq 0.
\end{align}
With the help of the global relation \eqref{gri}, we can modify the forcing so that the above equation now becomes:
\begin{align}\label{Delta2i2}
\Delta \mu_2&=-\frac{i}{2\pi}\int_{0}^{\infty}dl\int_{-\infty}^{\infty}d\xi\int_0^{y} d\eta\, e^{-il(\xi-x)-il(l+2k_R)(\eta-y)}q\Delta \mu_2\nonumber\\
&+\frac{1}{2\pi}\int_{0}^{\infty}dl\int_{-\infty}^{\infty}d\xi\int_{0}^{t}d\tau\,e^{-il(\xi-x)+il(l+2k_R)y+\omega(k_R,l)(\tau-t)}H(\xi,\tau,k_R,l)\Delta \phi_2\nonumber\\
&+\frac{i}{2\pi}\int_{-\infty}^{0}dl\int_{-\infty}^{\infty}d\xi\int_{y}^\infty d\eta\,e^{-il(\xi-x)-il(l+2k_R)(\eta-y)}q\Delta \mu_2\nonumber\\
&-\frac{i}{2\pi}\int_0^{\infty} dl\int_{-\infty}^{\infty}d\xi\int_{0}^{\infty}d\eta\,e^{-il(\xi-x)-il(l+2k_R)(\eta-y)-\omega(k_R,l)t}q_0\mu_2^-\Big |_{k_I=0,\, t=0}\nonumber\\
&+\frac{1}{2\pi}\int_{-\infty}^{-2k_R}dl\int_{-\infty}^{\infty}d\xi\int_{0}^{T}d\tau\,e^{-il(\xi-x)+il(l+2k_R)y+\omega(k_R,l)(\tau-t)}H(\xi,\tau,k_R,l)\phi_2^-\Big |_{k_I=0}\nonumber\\
&+\frac{i}{2\pi}\int_{-\infty}^0 dl\int_{-\infty}^{\infty}d\xi\int_{0}^{\infty}d\eta\,e^{-il(\xi-x)-il(l+2k_R)(\eta-y)-\omega(k_R,l)t}q_0\mu_2^-\Big |_{k_I=0,\, t=0},\quad  k_R\geq 0.
\end{align}
Furthermore, in the case of $k_R\geq 0$, employing the global relation \eqref{gri} in equations \eqref{plus1itrans} and \eqref{plus2itrans} and rearranging, we find the following expressions for $\breve \mu_1^+$ and $\breve \mu_2^+$:
\begin{align}\label{plus1i2trans}
&e_{2_{Xt}}\breve\mu_1^+=e_{2_{Xt}}-\frac{i}{2\pi}\int_{0}^{\infty}dl\int_{-\infty}^{\infty}d\xi\int_{0}^{y}d\eta\, e^{-il(\xi-x)-il(l+2k_R)(\eta-y)}q\left(e_{2_{\Xi t}}\breve\mu_1^+\right)\nonumber\\
&+\frac{1}{2\pi}\int_{0}^\infty dl\int_{-\infty}^{\infty}d\xi\int_{0}^t d\tau\,e^{-il(\xi-x)+il(l+2k_R)y+\omega(k_R,l)(\tau-t)}H(\xi,\tau,k_R,l)\left(e_{2_{\Xi \tau}}|_{\eta=0}\breve\phi_1^+\right)\nonumber\\
&+\frac{i}{2\pi}\int_{-\infty}^{0}dl\int_{-\infty}^{\infty}d\xi\int_{y}^{\infty}d\eta\,e^{-il(\xi-x)-il(l+2k_R)(\eta-y)}q\left(e_{2_{\Xi t}}\breve\mu_1^+\right)\nonumber\\
&+\int_{0}^{-2k_R-\lambda}dl\,e^{ilx+il(l+2k_R)y-\omega(k_R,l)t}\left[p_1^+(-k_R-\lambda,l+2k_R+\lambda)+r_1^+(-k_R-\lambda,l+2k_R+\lambda)\right],\nonumber\\
&\qquad \qquad \qquad \qquad \qquad \qquad \qquad \qquad \qquad \qquad \qquad k_R\geq 0, \,\lambda\geq -2k_R,
\end{align}
and
\begin{align}\label{plus2i2trans}
e_{2_{Xt}}\breve\mu_2^+&=e_{2_{Xt}}-\frac{i}{2\pi}\int_{0}^{\infty}dl\int_{-\infty}^{\infty}d\xi\int_{0}^{y}d\eta\, e^{-il(\xi-x)-il(l+2k_R)(\eta-y)}q\left(e_{2_{\Xi t}}\breve\mu_2^+\right)\nonumber\\
&+\frac{1}{2\pi}\int_{0}^\infty dl\int_{-\infty}^{\infty}d\xi\int_{0}^t d\tau\,e^{-il(\xi-x)+il(l+2k_R)y+\omega(k_R,l)(\tau-t)}H(\xi,\tau,k_R,l)\left(e_{2_{\Xi \tau}}|_{\eta=0}\breve\phi_2^+\right)\nonumber\\
&+\frac{i}{2\pi}\int_{-\infty}^{0}dl\int_{-\infty}^{\infty}d\xi\int_{y}^{\infty}d\eta\,e^{-il(\xi-x)-il(l+2k_R)(\eta-y)}q\left(e_{2_{\Xi t}}\breve\mu_2^+\right)\nonumber\\
&+\frac{i}{2\pi}\int_{0}^{-2k_R-\lambda}dl\,e^{ilx+il(l+2k_R)y-\omega(k_R,l)t}r_2^+(-k_R-\lambda,l+2k_R+\lambda), \quad k_R\geq 0, \,\lambda\leq -2k_R.
\end{align}
where the functions $p_1^+$, $r_1^+$ and $r_2^+$ are defined by equations \eqref{pj} and \eqref{rj}. 
\vskip 2mm
The above equations imply the following result:

\begin{prop}\label{D2}
The solution $\Delta\mu_2$ of equation \eqref{Delta2i} is given by
\begin{align}\label{solD2}
\Delta\mu_2(x,y,t,k_R)&=\int_{-\infty}^{-2k_R}d\lambda\, \psi_1(k_R,\lambda) \left(e_{2_{Xt}}\breve \mu_2^+\right)+\int_{-2k_R}^{\infty}d\lambda\, \psi_2(k_R,\lambda) \left(e_{2_{Xt}}\breve \mu_1^+\right)\nonumber\\
&+\int_{-2k_R}^0d\lambda\, \psi_3(k_R,\lambda) \left(e_{2_{Xt}}\breve \mu_2^+\right), \quad k_R\geq 0,,
\end{align}
provided that the functions $\psi_j,\ j=1,2,3,$ satisfy the following linear Volterra integral equations: 
\begin{subequations}\label{psi}
\begin{align}
&\psi_1(k_R,\lambda)+\int_{-\infty}^{\lambda}dl\, \psi_1(k_R,l)r_2^+(-k_R-l,l-\lambda)=-r_2^-(k_R,-\lambda-2k_R), \quad k_R\geq 0,\ \lambda \leq -2k_R,\label{psi1}\\
&\psi_2(k_R,\lambda)-\int_{\lambda}^\infty dl\, \psi_2(k_R,l)\left[p_1^+(-k_R-l,l-\lambda)+r_1^+(-k_R-l,l-\lambda)\right]=\nonumber\\
&=p_2^-(k_R,l-\lambda)+r_2^-(k_R,-\lambda-2k_R), \quad k_R\geq 0,\ \lambda \geq -2k_R,\label{psi2}\\ 
\nonumber\\
&\psi_3(k_R,\lambda)-\int_{\lambda}^{0}dl\, \psi_3(k_R,l)\left[p_1^+(-k_R-l,l-\lambda)+r_1^+(-k_R-l,l-\lambda)\right]=\nonumber\\
&=-p_2^-(k_R,-\lambda-2k_R), \quad k_R\geq 0,\ -2k_R\leq\lambda \leq 0,\label{psi3}.  
\end{align}
\end{subequations}
where the functions $p_j^\pm$ and $r_j^\pm$, $j=1,2,$ are defined by equations \eqref{pj} and \eqref{rj} respectively.
\end{prop}

\begin{PROOF}
The forcing terms of $\mathcal S_2$ are equal to
\begin{align}\label{fS2}
&\int_{-\infty}^{-2k_R}d\lambda\, \psi_1(k_R,\lambda)\,e_{2_{Xt}}+\int_{-\infty}^{-2k_R}d\lambda\, \psi_1(k_R,\lambda)\int_{0}^{-2k_R-\lambda} dl\,E_{Xt}r_2^+(-k_R-\lambda,l+2k_R+\lambda)\nonumber\\
&+\int_{-2k_R}^\infty d\lambda\, \psi_2(k_R,\lambda)e_{2_{Xt}}\nonumber\\
&+\int_{-2k_R}^\infty d\lambda\, \psi_2(k_R,\lambda)\int_{0}^{-2k_R-\lambda} dl\,E_{Xt}\left[p_1^+(-k_R-\lambda,l+2k_R+\lambda)+r_1^+(-k_R-\lambda,l+2k_R+\lambda)\right]\nonumber\\
&+\int^{0}_{-2k_R}d\lambda\, \psi_3(k_R,\lambda)\,e_{2_{Xt}}\nonumber\\
&+\int^{0}_{-2k_R}d\lambda\, \psi_3(k_R,\lambda)\int^{0}_{-2k_R} dl\,E_{Xt}\left[p_1^+(-k_R-\lambda,l+2k_R+\lambda)+r_1^+(-k_R-\lambda,l+2k_R+\lambda)\right],
\end{align}
where $e_{2_{Xt}}$ and $E_{Xt}$ are defined by equations \eqref{eXt} and \eqref{EXt} respectively. Moreover, interchanging the order of integration between $\lambda$ and $l$, the forcing becomes
\begin{align}\label{fS22}
&\int_{-\infty}^{-2k_R}d\lambda\, \psi_1(k_R,\lambda)\,e_{2_{Xt}}+\int_{0}^{\infty}dl\,E_{Xt}\int_{0}^{-2k_R-l} d\lambda\, \psi_1(k_R,\lambda)r_2^+(-k_R-\lambda,l+2k_R+\lambda)\nonumber\\
&+\int_{-2k_R}^\infty d\lambda\, \psi_2(k_R,\lambda)e_{2_{Xt}}\nonumber\\
&-\int_{-\infty}^0 dl E_{Xt}\int_{-2k_R-l}^{\infty} d\lambda\,\psi_2(k_R,\lambda)\left[p_1^+(-k_R-\lambda,l+2k_R+\lambda)+r_1^+(-k_R-\lambda,l+2k_R+\lambda)\right]\nonumber\\
&+\int^{0}_{-2k_R}d\lambda\, \psi_3(k_R,\lambda)\,e_{2_{Xt}}\nonumber\\
&-\int^{0}_{-2k_R}dl\,E_{Xt} \int^{-2k_R-l}_{0} d\lambda\,\psi_3(k_R,\lambda)\left[p_1^+(-k_R-\lambda,l+2k_R+\lambda)+r_1^+(-k_R-\lambda,l+2k_R+\lambda)\right],\quad k_R\geq 0.
\end{align}
Under the change of variables $l\mapsto -2k_R-l$ in the integrals involving $E_{Xt}$, equation \eqref{fS22} takes the following form:
\begin{align}\label{fS23}
&\int_{-\infty}^{-2k_R}dl\, \psi_1(k_R,l)\,e_{2_{Xt}}+\int_{-\infty}^{-2k_R}dl\,e_{2_{Xt}}\int_{-\infty}^{l} d\lambda\, \psi_1(k_R,\lambda)r_2^+(-k_R-\lambda,\lambda-l)\nonumber\\
&+\int_{-2k_R}^\infty dl\, \psi_2(k_R,l)e_{2_{Xt}}\nonumber\\
&-\int_{-2k_R}^\infty dl e_{2_{Xt}}\int_{l}^{\infty} d\lambda\,\psi_2(k_R,\lambda)\left[p_1^+(-k_R-\lambda,\lambda-l)+r_1^+(-k_R-\lambda,\lambda-l)\right]\nonumber\\
&+\int^{0}_{-2k_R}dl\, \psi_3(k_R,l)\,e_{2_{Xt}}\nonumber\\
&-\int^{0}_{-2k_R}dl\,e_{2_{Xt}} \int^{l}_{0} d\lambda\,\psi_3(k_R,\lambda)\left[p_1^+(-k_R-\lambda,\lambda-l)+r_1^+(-k_R-\lambda,\lambda-l)\right],\quad k_R\geq 0.
\end{align}
Furthermore, under the change of variables $l\mapsto -2k_R-l$ and with the aid of the definitions \eqref{pjrj}, the forcing of equation \eqref{Delta2i2} is equal to:
\begin{align}\label{fD2}
&-\int^{0}_{-2k_R}dl\, e_{2_{Xt}}\, p_2^-(k_R,-2k_R-l)-\int_{-\infty}^{-2k_R}dl\, e_{2_{Xt}}\, r_2^-(k_R,-2k_R-l)\nonumber\\
&+\int_{-2k_R}^{\infty}dl\, e_{2_{Xt}}\left[p_2^-(k_R,-2k_R-l)+r_2^-(k_R,-2k_R-l)\right].
\end{align}
For equations \eqref{Delta2i2} and \eqref{solD2} to be identical, their forcing terms \eqref{fS23} and \eqref{fD2} must be equal, therefore, we define the functions $\psi_j,\, j=1,2,3,$ via equations \eqref{psi}. In this case, uniqueness implies that $\Delta \mu_2=\mathcal S_2$. 
\QED
\end{PROOF}

\subsection{The spectral functions}\label{kpimap}

We will now define the map 
\begin{align}\label{mapi}
\{q_0(x,y),\ g(x,t),\ h(x,t)\}\mapsto\{f_1^+(k_R,k_I),f_2^- (k_R,k_I), \chi_{j}(k_R,\lambda),\psi_{j}(k_R,\lambda)\},\quad j=1,2,3,
\end{align}
from the initial condition $q_0$ and the boundary values $g$ and $h$ to the spectral functions that are defined by equations \eqref{fplus1i}, \eqref{fminus2i}, \eqref{chi} and \eqref{psi}.

\begin{enumerate}

\item {\small$\{q_0,g,h\} \mapsto \{\rho_{j}^\pm(x,y,k_R,k_I),\ \phi_{j}^\pm(x,t,k_R,k_I)\}, \quad j=1,2$.}

\noindent The functions $\rho_{j}^\pm$ are defined in terms of $\{q_0,g,h\}$ and $\phi_{j}^\pm$ via the following linear integral equations:
\begin{subequations}\label{maprhoi}
\begin{align}\label{rhoplus1i}
\rho_1^+&=1-\frac{i}{2\pi}\int_{0}^{\infty}dl\int_{-\infty}^{\infty}d\xi\int_{0}^{y}d\eta\, e^{-il(\xi-x)-il(l+2k)(\eta-y)}q_0\rho_1^+\nonumber\\
&-\frac{1}{2\pi}\int_{0}^{-2k_R}dl\int_{-\infty}^{\infty}d\xi\int_{0}^Td\tau\,e^{-il(\xi-x)+il(l+2k)y+\omega(k,l)\tau}H\phi_1^+\nonumber\\
&+\frac{i}{2\pi}\int_{-\infty}^{0}dl\int_{-\infty}^{\infty}d\xi\int_{y}^{\infty}d\eta\,e^{-il(\xi-x)-il(l+2k)(\eta-y)}q_0\rho_1^+, \quad k_R\leq 0,\, k_I\geq 0,
\end{align}
\begin{align}\label{rhominus1i}
\rho_1^-&=1+\frac{i}{2\pi}\int_{0}^{\infty}dl\int_{-\infty}^{\infty}d\xi\int_{y}^\infty d\eta\, e^{-il(\xi-x)-il(l+2k)(\eta-y)}q_0\rho_1^-\nonumber\\
&-\frac{1}{2\pi}\int_{-\infty}^0dl\int_{-\infty}^{\infty}d\xi\int_{0}^{T}d\tau\,e^{-il(\xi-x)+il(l+2k)y+\omega(k,l)\tau}H\phi_1^-\nonumber\\
&-\frac{i}{2\pi}\int_{-\infty}^{0}dl\int_{-\infty}^{\infty}d\xi\int_{0}^{y}d\eta\,e^{-il(\xi-x)-il(l+2k)(\eta-y)}q_0\rho_1^-,\quad  k_R\leq 0,\, k_I\leq 0,
\end{align}
\begin{align}\label{rhominus2i}
\rho_2^-&=1+\frac{i}{2\pi}\int_{0}^{\infty}dl\int_{-\infty}^{\infty}d\xi\int_{y}^\infty d\eta\, e^{-il(\xi-x)-il(l+2k)(\eta-y)}q_0\rho_2^-\nonumber\\
&-\frac{1}{2\pi}\int_{-\infty}^0dl\int_{-\infty}^{\infty}d\xi\int_{0}^{T}d\tau\,e^{-il(\xi-x)+il(l+2k)y+\omega(k,l)\tau}H\phi_2^-\nonumber\\
&-\frac{i}{2\pi}\int_{-\infty}^{0}dl\int_{-\infty}^{\infty}d\xi\int_{0}^{y}d\eta\,e^{-il(\xi-x)-il(l+2k)(\eta-y)}q_0\rho_2^-,\quad  k_R\geq 0,\, k_I\leq 0,
\end{align}
\begin{align}\label{rhoplus2i}
\rho_2^+&=1-\frac{i}{2\pi}\int_{0}^{\infty}dl\int_{-\infty}^{\infty}d\xi\int_{0}^{y}d\eta\, e^{-il(\xi-x)-il(l+2k)(\eta-y)}q_0\rho_2^+\nonumber\\
&+\frac{i}{2\pi}\int_{-\infty}^{0}dl\int_{-\infty}^{\infty}d\xi\int_{y}^{\infty}d\eta\,e^{-il(\xi-x)-il(l+2k)(\eta-y)}q_0\rho_2^+, \quad k_R\geq 0,\, k_I\geq 0.
\end{align}
\end{subequations}
Evaluating equations \eqref{plus1i}-\eqref{plus2i} at $y=0$, we obtain the following integral equations for $\phi_{j}^\pm$ in terms of $q$, $g$ and $h$:
\begin{subequations}\label{phii}
\begin{align}\label{phplus1i}
\phi_1^+&=1-\frac{1}{2\pi}\int_{0}^{-2k_R}dl\int_{-\infty}^{\infty}d\xi\int_{t}^Td\tau\,e^{-il(\xi-x)+\omega(k,l)(\tau-t)}H\phi_1^+\nonumber\\
&+\frac{1}{2\pi}\int_{-2k_R}^{\infty}dl\int_{-\infty}^{\infty}d\xi\int_{0}^{t}d\tau\,e^{-il(\xi-x)+\omega(k,l)(\tau-t)}H\phi_1^+\nonumber\\
&+\frac{i}{2\pi}\int_{-\infty}^{0}dl\int_{-\infty}^{\infty}d\xi\int_{0}^{\infty}d\eta\,e^{-il(\xi-x)-il(l+2k)\eta}q\mu_1^+, \quad k_R\leq 0,\, k_I\geq 0,
\end{align}
\begin{align}\label{phminus1i}
\phi_1^-&=1+\frac{i}{2\pi}\int_{0}^{\infty}dl\int_{-\infty}^{\infty}d\xi\int_{0}^\infty d\eta\, e^{-il(\xi-x)-il(l+2k)\eta}q\mu_1^-\nonumber\\
&-\frac{1}{2\pi}\int_{-\infty}^0dl\int_{-\infty}^{\infty}d\xi\int_{t}^{T}d\tau\,e^{-il(\xi-x)+\omega(k,l)(\tau-t)}H\phi_1^-,\quad  k_R\leq 0,\, k_I\leq 0,
\end{align}
\begin{align}\label{phminus2i}
\phi_2^-&=1+\frac{i}{2\pi}\int_{0}^{\infty}dl\int_{-\infty}^{\infty}d\xi\int_{0}^\infty d\eta\, e^{-il(\xi-x)-il(l+2k)\eta}q\mu_2^-\nonumber\\
&-\frac{1}{2\pi}\int_{-\infty}^0dl\int_{-\infty}^{\infty}d\xi\int_{t}^{T}d\tau\,e^{-il(\xi-x)+\omega(k,l)(\tau-t)}H\phi_2^-\nonumber\\
&+\frac{1}{2\pi}\int_{-2k_R}^0dl\int_{-\infty}^{\infty}d\xi\int_{0}^{t}d\tau\,e^{-il(\xi-x)+\omega(k,l)(\tau-t)}H\phi_2^-,\quad  k_R\geq 0,\, k_I\leq 0,
\end{align}
and
\begin{align}\label{phplus2i}
\phi_2^+&=1+\frac{1}{2\pi}\int_{0}^{\infty}dl\int_{-\infty}^{\infty}d\xi\int_{0}^{t}d\tau\,e^{-il(\xi-x)+\omega(k,l)(\tau-t)}H\phi_2^+\nonumber\\
&+\frac{i}{2\pi}\int_{-\infty}^{0}dl\int_{-\infty}^{\infty}d\xi\int_{0}^{\infty}d\eta\,e^{-il(\xi-x)-il(l+2k)\eta}q\mu_2^+, \quad k_R\geq 0,\, k_I\geq 0,
\end{align}
\end{subequations}
with $H$ defined by equation \eqref{H}.

We must eliminate $q$ and $\mu_j^\pm$ from the above expressions. Since the global relation \eqref{gri}, is valid only for $\mathrm {Re}\, \omega(k,l)\geq 0$, we will employ its alternative form \eqref{gr2i}, which is valid in the case of $\mathrm {Re}\, \omega(k,l)\leq 0$, i.e. for $lk_I\leq0$ and $lk_I(l+2k_R)\leq 0$:

\begin{align*}
&\int_{-\infty}^{\infty}d\xi\int_{0}^\infty d\eta\, e^{-il\xi-il(l+2k)\eta} q(\xi,\eta,t)\mu(\xi,\eta,t,k_R,k_I)\nonumber\\
&+i\int_{-\infty}^{\infty}d\xi\int_{t}^{T}d\tau\, e^{-il\xi+\omega(k,l)(\tau-t)}H(\xi,\tau,k,l)\phi(\xi,\tau,k_R,k_I)\nonumber\\
&=\int_{-\infty}^{\infty}d\xi\int_{0}^\infty d\eta\, e^{-il\xi-il(l+2k)\eta+\omega(k,l)(T-t)} q(\xi,\eta,T)\mu(\xi,\eta,T,k_R,k_I).
\end{align*}
The unknown $q(x,y,T)$ can be eliminated by taking the limit $T\rightarrow \infty$, so that the RHS of the above equation vanishes due to exponential decay. 

Hence, the functions $\phi^\pm_j$ can be defined in terms of $q_0$, $g$, $h$ and $\rho_j^\pm$ via the following integral equations: 
\begin{subequations}\label{phiimap}
\begin{align}\label{phiplus1i}
\phi_1^+&=1-\frac{1}{2\pi}\int_{0}^{-2k_R}dl\int_{-\infty}^{\infty}d\xi\int_{t}^\infty d\tau\,e^{-il(\xi-x)+\omega(k,l)(\tau-t)}H\phi_1^+\nonumber\\
&+\frac{1}{2\pi}\int_{-2k_R}^{\infty}dl\int_{-\infty}^{\infty}d\xi\int_{0}^{t}d\tau\,e^{-il(\xi-x)+\omega(k,l)(\tau-t)}H\phi_1^+\nonumber\\
&-\frac{1}{2\pi}\int_{-\infty}^{0}dl\int_{-\infty}^{\infty}d\xi\int_{0}^{t}d\tau\,e^{-il(\xi-x)+\omega(k,l)(\tau-t)}H\phi_1^+\nonumber\\
&+\frac{i}{2\pi}\int_{-\infty}^{0}dl\int_{-\infty}^{\infty}d\xi\int_{0}^\infty d\eta\, e^{-il(\xi-x)-il(l+2k)\eta-\omega(k,l)t} q_0\rho_1^+, \quad k_R\leq 0,\, k_I\geq 0,
\end{align}
\begin{align}\label{phiminus1i}
\phi_1^-&=1+\frac{1}{2\pi}\int_{0}^{-2k_R}dl\int_{-\infty}^{\infty}d\xi\int_{0}^t d\tau\, e^{-il(\xi-x)+\omega(\tau-t)} H\phi_1^+\nonumber\\
&+\frac{1}{2\pi}\int_{-2k_R}^\infty dl\int_{-\infty}^{\infty}d\xi\int_{t}^\infty d\tau\, e^{-il(\xi-x)+\omega(\tau-t)} H\phi_1^+\nonumber\\
&-\frac{1}{2\pi}\int_{-\infty}^0dl\int_{-\infty}^{\infty}d\xi\int_{t}^{\infty}d\tau\,e^{-il(\xi-x)+\omega(k,l)(\tau-t)}H\phi_1^-\nonumber\\
&+\frac{i}{2\pi}\int_{0}^{-2k_R}dl\int_{-\infty}^{\infty}d\xi\int_{0}^\infty d\eta\, e^{-il(\xi-x)-il(l+2k)\eta-\omega(k,l)t}q_0\rho_1^-,\quad  k_R\leq 0,\, k_I\leq 0,
\end{align}
\begin{align}\label{phiminus2i}
\phi_2^-&=1+\frac{1}{2\pi}\int_{0}^{\infty}dl\int_{-\infty}^{\infty}d\xi\int_{t}^\infty d\tau\, e^{-il(\xi-x)+\omega(k,l)(\tau-t)}H\phi_2^-\nonumber\\
&-\frac{1}{2\pi}\int_{-\infty}^0dl\int_{-\infty}^{\infty}d\xi\int_{t}^{\infty}d\tau\,e^{-il(\xi-x)+\omega(k,l)(\tau-t)}H\phi_2^-\nonumber\\
&+\frac{1}{2\pi}\int_{-2k_R}^0dl\int_{-\infty}^{\infty}d\xi\int_{0}^{t}d\tau\,e^{-il(\xi-x)+\omega(k,l)(\tau-t)}H\phi_2^-,\quad  k_R\geq 0,\, k_I\leq 0,
\end{align}
and
\begin{align}\label{phiplus2i}
\phi_2^+&=1+\frac{1}{2\pi}\int_{0}^{\infty}dl\int_{-\infty}^{\infty}d\xi\int_{0}^{t}d\tau\,e^{-il(\xi-x)+\omega(k,l)(\tau-t)}H\phi_2^+\nonumber\\
&-\frac{1}{2\pi}\int_{-\infty}^{-2k_R}dl\int_{-\infty}^{\infty}d\xi\int_{0}^{t}d\tau\,e^{-il(\xi-x)+\omega(k,l)(\tau-t)}H\phi_2^+\nonumber\\
&+\frac{1}{2\pi}\int_{-2k_R}^{0}dl\int_{-\infty}^{\infty}d\xi\int_{t}^\infty d\tau\,e^{-il(\xi-x)+\omega(k,l)(\tau-t)}H\phi_2^+\nonumber\\
&+\frac{i}{2\pi}\int_{-\infty}^0 dl\int_{-\infty}^{\infty}d\xi\int_{0}^{\infty}d\eta\,e^{-il(\xi-x)-il(l+2k)\eta}q_0\rho_2^+, \quad k_R\geq 0,\, k_I\geq 0.
\end{align}
\end{subequations}

\item{$\{q_0,g,h,\phi_1^+, \phi_2^-\}\mapsto \{f_1^+, f_2^-\}$,}

where the functions $f_1^+$ and $f_2^-$ are defined by:
\begin{subequations}\label{mapf1f2}
\begin{align}
f_1^+(k_R,k_I)&=\frac{1}{2\pi}\int_{-\infty}^{\infty}d\xi\int_{0}^{T}d\tau\,e^{2ik_R\xi+\omega(k,-2k_R)\tau}H(\xi,\tau,k,-2k_R)\phi_1^+
\label{mapfplus1i},\\
f_2^-(k_R,k_I)&=\frac{1}{2\pi}\int_{-\infty}^{\infty}d\xi\int_{0}^{T}d\tau\,e^{2ik_R\xi+\omega(k,-2k_R)\tau}H(\xi,\tau,k,-2k_R)\phi_2^-.\label{mapfminus2i}
\end{align}
\end{subequations}

\item{$\{q_0,g,h,\phi_{j}^\pm, \rho_{j}^\pm\} \mapsto \{p_{j}^\pm,\ r_{j}^\pm \}, \quad j=1,2,$} 

where the functions $p_{j}^\pm$ and $r_j^\pm$ are given by:
\begin{subequations}\label{mapr1r2p1p2}
\begin{align}
p_j^\pm(k,l)&=\frac{1}{2\pi}\int_{-\infty}^{\infty}d\xi\int_{0}^{T}d\tau\,e^{-il\xi+\omega(k,l)\tau}H \phi_j^\pm, &\quad j=1,2,\label{mapp12}\\
r_j^\pm(k,l)&=\frac{i}{2\pi}\int_{-\infty}^{\infty}d\xi\int_{0}^{\infty}d\eta\,e^{-il\xi-il(l+2k)\eta}q_0\mu_j^\pm\Big|_{t=0}, &\quad j=1,2.\label{mapr12}
\end{align}
\end{subequations}

\item{$\{p_j^\pm, r_j^\pm, \ j=1,2\}\mapsto \{\chi_j,\psi_j, \ j=1,2,3\}$,}

where the functions $\chi_{j}$ and $\psi_{j}$ are defined via the following Volterra integral equations:
\begin{subequations}\label{mapchi}
\begin{align}
&\chi_1(k_R,\lambda)+\int_{-\infty}^{\lambda}dl\, \chi_1(k_R,l)r_2^+(-k_R-l,l-\lambda)=-r_1^-(k_R,-\lambda-2k_R), \quad k_R\leq 0,\ \lambda \leq -2k_R,\label{mapchi1}\\
&\chi_2(k_R,\lambda)-\int_{\lambda}^\infty dl\, \chi_2(k_R,l)\left[p_1^+(-k_R-l,l-\lambda)+r_1^+(-k_R-l,l-\lambda)\right]=\nonumber\\
&=p_1^-(k_R,-\lambda-2k_R)+r_1^-(k_R,-\lambda-2k_R), \quad k_R\leq 0,\ \lambda \geq -2k_R,\label{mapchi2}\\ 
\nonumber\\
&\chi_3(k_R,\lambda)+\int_{0}^{-2k_R}dl\, \chi_3(k_R,l)p_2^+(-k_R-l,l-\lambda)\nonumber\\
&+\int_{0}^{\lambda}dl\, \chi_3(k_R,l)r_2^+(-k_R-l,l-\lambda)+\int_{-\infty}^{-2k_R}dl\, \chi_1(k_R,l)p_2^+(-k_R-l,l-\lambda)\nonumber\\
&+\int_{-2k_R}^{\infty}dl\, \chi_2(k_R,l)p_1^+(-k_R-l,l-\lambda)=-p_1^-(k_R,-\lambda-2k_R), \quad k_R\leq 0,\ 0\leq\lambda \leq -2k_R,\label{mapchi3}
\end{align}
\end{subequations}
and
\begin{subequations}\label{mappsi}
\begin{align}
&\psi_1(k_R,\lambda)+\int_{-\infty}^{\lambda}dl\, \psi_1(k_R,l)r_2^+(-k_R-l,l-\lambda)=-r_2^-(k_R,-\lambda-2k_R), \quad k_R\geq 0,\ \lambda \leq -2k_R,\label{mappsi1}\\
&\psi_2(k_R,\lambda)-\int_{\lambda}^\infty dl\, \psi_2(k_R,l)\left[p_1^+(-k_R-l,l-\lambda)+r_1^+(-k_R-l,l-\lambda)\right]=\nonumber\\
&=p_2^-(k_R,-\lambda-2k_R)+r_2^-(k_R,-\lambda-2k_R), \quad k_R\geq 0,\ \lambda \geq -2k_R,\label{mappsi2}\\ 
\nonumber\\
&\psi_3(k_R,\lambda)-\int_{\lambda}^{0}dl\, \psi_3(k_R,l)\left[p_1^+(-k_R-l,l-\lambda)+r_1^+(-k_R-l,l-\lambda)\right]=\nonumber\\
&=-p_2^-(k_R,-\lambda-2k_R), \quad k_R\geq 0,\ -2k_R\leq\lambda \leq 0\label{mappsi3}.  
\end{align}
\end{subequations}
\end{enumerate}
\paragraph{Remark 4.3.} Note that, since we need to take $T\rightarrow \infty$ in order to define the functions $\rho_j, \phi_j$, $j=1,2,$ of mapping 1, every function involved in the mappings 2-4 should be defined under the limit $T\rightarrow \infty$. 

\begin{prop}
The solution $q(x,y,t)$ to the initial-boundary value problem for the KPI equation \eqref{kpiintro} admits the following integral representation:

\begin{align}\label{qprop}
\frac{i\pi}{2}\,q&=\int_{-\infty}^0\!\!d\nu_R \int_{-i \infty}^{i\nu_R}\!d\nu_I\,\psi_1(\nu,-2i \nu_I)\Bigg(e_{1_{Xt}}(x,y,t,\nu_R,\nu_I)\mu_2^+(x,y,t,- \nu_R, \nu_I)\Bigg)_x\nonumber\\
&+ \int_{0}^{\infty} d\nu_R\int^{i\nu_R}_{-i\infty}d\nu_I\,\psi_2(\nu,-2i\nu_I)\Bigg(e_{1_{Xt}}(x,y,t,\nu_R,\nu_I) \mu_1^+(x,y,t,-\nu_R,\nu_I)\Bigg)_x\nonumber\\
&+ \int_{0}^{\infty} d\nu_R\int_{-i\infty}^0 d\nu_I\,\psi_3(\nu,-2i\nu_I)\Bigg(e_{1_{Xt}}(x,y,t,\nu_R,\nu_I)\mu_2^+(x,y,t,-\nu_R,\nu_I)\Bigg)_x\nonumber\\
&+ \int_{-\infty}^0d \nu_R \int_{i \nu_R}^{i\infty}d \nu_I\,\chi_1( \nu,-2i \nu_I)\Bigg(e_{1_{Xt}}(x,y,t,\nu_R,\nu_I)\mu_1^+(x,y,t,- \nu_R, \nu_I)\Bigg)_x\nonumber\\
&+ \int_{-\infty}^0 d \nu_R\int_{i \nu_R}^{i \infty}d\nu_I\,\chi_2( \nu,-2i \nu_I)\Bigg(e_{1_{Xt}}(x,y,t,\nu_R,\nu_I)\mu_1^+(x,y,t,- \nu_R, \nu_I)\Bigg)_x\nonumber\\
&+\int_{-\infty}^0d\nu_R\int_{i \nu_R}^{0}d\nu_I\,\chi_3( \nu,-2i \nu_I)\Bigg(e_{1_{Xt}}(x,y,t,\nu_R,\nu_I)\mu_2^+(x,y,t,- \nu_R, \nu_I)\Bigg)_x\nonumber\\
&+\int_{-\infty}^0 d\nu_R\int_{0}^\infty d\nu_I\,f_1^{+}(\nu_R,\nu_I,T)\Bigg(e_{1_{Xt}}(x,y,t,\nu_R,\nu_I) \mu_2^+(x,y,t,-\nu_R, \nu_I)\Bigg)_x\nonumber\\
&+\int_{0}^\infty d\nu_R\int_{-\infty}^{0} d\nu_I\,f_2^{-}(\nu_R,\nu_I,T)\Bigg(e_{1_{Xt}}(x,y,t,\nu_R,\nu_I) \mu_1^-(x,y,t,-\nu_R, \nu_I)\Bigg)_x,
\end{align}
where
\begin{equation*}
e_{1_{Xt}}(x,y,t,k_R,k_I)=e^{-2ik_R x+4k_Rk_I y-\omega(k,-2k_R)t},
\end{equation*}
and the functions $f_1^+$, $f_2^-$, $\chi_j$ and $\psi_j$, $j=1,2,3,$ are defined via equations \eqref{f}, \eqref{mapchi} and \eqref{mappsi}. 
\end{prop}

\begin{PROOF}
Substituting in the Pompeiu's formula \eqref{pompi} the expressions for the d-bar derivatives $\partial \mu^\pm_j/\partial \bar k,\ j=1,2,$ (see proposition \ref{kpid-barprop}) and for the jumps $\Delta \mu_j,\ j=1,2,$ (see propositions \ref{D1} and \ref{D2}), we obtain the following representation for $\mu$:
\begin{align}\label{pompi2}
\mu&=1+\frac{1}{2i \pi} \int_{0}^{\infty} \frac{d\nu_R}{\nu_R-k} \left[\int_{-\infty}^{-2\nu_R}d\lambda\, \psi_1(\nu_R,\lambda) e_{2_{Xt}} \breve \mu_2^+\right.\nonumber\\
&\left.+\int_{-2\nu_R}^{\infty} d\lambda\, \psi_2(\nu_R,\lambda) e_{2_{Xt}} \breve \mu_1^++\int_{-2\nu_R}^0 d\lambda\, \psi_3(\nu_R,\lambda) e_{2_{Xt}} \breve \mu_2^+\right]\nonumber\\
&+\frac{1}{2i \pi} \int_{-\infty}^0 \frac{d\nu_R}{\nu_R-k} \left[\int_{-\infty}^{-2\nu_R}d\lambda\, \chi_1(\nu_R,\lambda) e_{2_{Xt}} \breve \mu_2^+\right.\nonumber\\
&\left.+\int_{-2\nu_R}^{\infty} d\lambda\, \chi_2(\nu_R,\lambda) e_{2_{Xt}} \breve \mu_1^++\int_{0}^{-2\nu_R} d\lambda\, \chi_3(\nu_R,\lambda) e_{2_{Xt}} \breve \mu_2^+\right]\nonumber\\
&-\frac{1}{\pi}\int_{-\infty}^{0}d\nu_R\int_0^\infty\!\! \frac{d\nu_I}{\nu-k}\,f_1^+(\nu_R,\nu_I,T)e_{1_{Xt}}\tilde\mu_2^+\nonumber\\
&-\frac{1}{\pi}\int_0^{\infty}d\nu_R\int_{-\infty}^0 \frac{d\nu_I}{\nu-k}\,f_2^-(\nu_R,\nu_I,T)e_{1_{Xt}}\tilde\mu_1^-,
\end{align}
where $\tilde \mu_j$, $e_{2_{Xt}}$ and $\breve \mu_j$ are defined by equations \eqref{mutilde}, \eqref{eXt2} and \eqref{muplus1breve}.

The change of variables
\begin{equation}\label{changevar}
\nu_R+\frac{\lambda}{2}=\hat \nu_R, \ \ \frac{\lambda}{2}=\hat \nu_I, 
\end{equation}
in the integrals with respect to $(\nu_R,\lambda)$ on the RHS of the above expression implies that 
\begin{equation}
 e_{2_{Xt}}(x,y,t,\nu_R,\lambda)\mapsto e_{Xt}(x,y,t,\hat \nu_R-\hat \nu_I,2\hat \nu_I):=e_{3_{Xt}},
\end{equation}
where $e_{Xt}$ is defined by \eqref{eXtgen}, hence (see also figures \ref{change} and \ref{change2}),

\begin{figure}[ht]
\begin{center}
\resizebox{8cm}{!}{\input{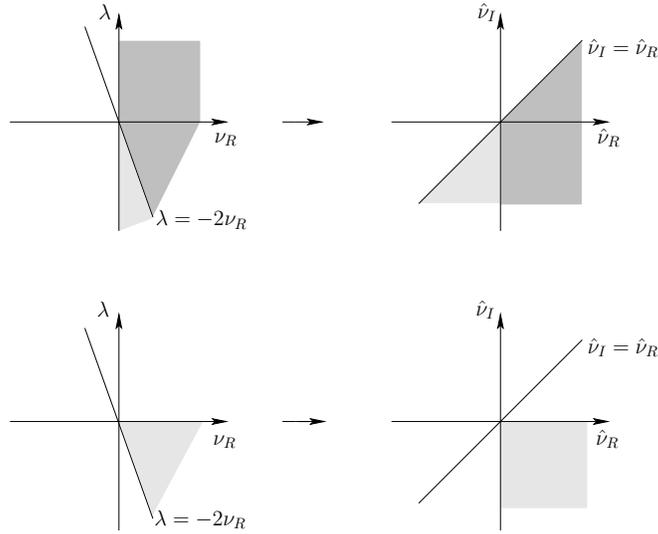}}
\end{center}
\caption{The change of variables \eqref{changevar} for $\psi_j$, $j=1,2,3$.}
\label{change}
\end{figure} 

\begin{figure}[ht]
\begin{center}
\resizebox{8cm}{!}{\input{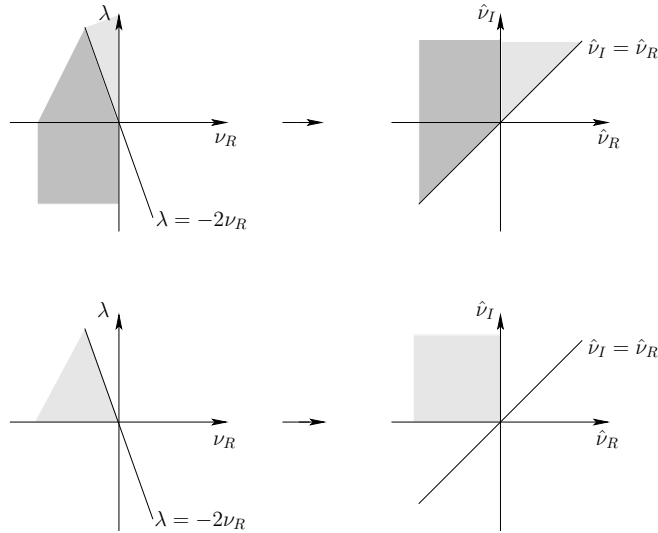}}
\end{center}
\caption{The change of variables \eqref{changevar} for $\chi_j$, $j=1,2,3$.}
\label{change2}
\end{figure}

\begin{align}\label{pompi3}
\mu&=1+\frac{1}{i\pi} \int_{-\infty}^{0} d\hat \nu_R\int_{-\infty}^{\hat\nu_R}\frac{d\hat \nu_I}{\hat \nu_R-\hat \nu_I-k}\, \psi_1(\hat \nu_R-\hat \nu_I,2\hat \nu_I) e_{3_{Xt}}  \mu_2^+(x,y,t,-\hat \nu_R,i \hat \nu_I)\nonumber\\
&+\frac{1}{i\pi}\int_{0}^\infty d\hat \nu_R\int_{-\infty}^{\hat\nu_R}\frac{d\hat \nu_I}{\hat \nu_R-\hat \nu_I-k}\, \psi_2(\hat \nu_R-\hat \nu_I,2\hat \nu_I) e_{3_{Xt}}  \mu_1^+(x,y,t,-\hat \nu_R,i \hat \nu_I)\nonumber\\
&+\frac{1}{i\pi}\int_{0}^\infty d\hat \nu_R\int_{-\infty}^{0}\frac{d\hat \nu_I}{\hat \nu_R-\hat \nu_I-k}\, \psi_3(\hat \nu_R-\hat \nu_I,2\hat \nu_I) e_{3_{Xt}}  \mu_2^+(x,y,t,-\hat \nu_R,i \hat \nu_I)\nonumber\\
&+\frac{1}{i\pi}\int_{-\infty}^0 d\hat \nu_R\int_{\hat \nu_R}^\infty \frac{d\hat \nu_I}{\hat \nu_R-\hat \nu_I-k}\, \chi_1(\hat \nu_R-\hat \nu_I,2\hat \nu_I) e_{3_{Xt}}  \mu_2^+(x,y,t,-\hat \nu_R,i \hat \nu_I)\nonumber\\
&+\frac{1}{i\pi}\int_{0}^\infty d\hat \nu_R\int_{\hat \nu_R}^\infty \frac{d\hat \nu_I}{\hat \nu_R-\hat \nu_I-k}\, \chi_2(\hat \nu_R-\hat \nu_I,2\hat \nu_I) e_{3_{Xt}}  \mu_1^+(x,y,t,-\hat \nu_R,i \hat \nu_I)\nonumber\\
&+\frac{1}{i\pi}\int_{-\infty}^0 d\hat \nu_R\int_{\hat \nu_R}^{0}\frac{d\hat \nu_I}{\hat \nu_R-\hat \nu_I-k}\, \chi_3(\hat \nu_R-\hat \nu_I,2\hat \nu_I) e_{3_{Xt}}  \mu_2^+(x,y,t,-\hat \nu_R,i \hat \nu_I)\nonumber\\
&-\frac{1}{\pi}\int_{-\infty}^{0}d\nu_R\int_0^\infty\!\! \frac{d\nu_I}{\nu-k}\,f_1^+(\nu_R,\nu_I,T)e_{1_{Xt}}\tilde\mu_2^+\nonumber\\
&-\frac{1}{\pi}\int_0^{\infty}d\nu_R\int_{-\infty}^0 \frac{d\nu_I}{\nu-k}\,f_2^-(\nu_R,\nu_I,T)e_{1_{Xt}}\tilde\mu_1^-.
\end{align}
Furthermore, letting $\tilde\nu_I=i \hat \nu_I$ and then dropping the hat and the tilde, we have:
{\small\begin{align}\label{pompi4}
\mu&=1-\frac{1}{\pi} \int_{-\infty}^{0} d\nu_R\int_{-i\infty}^{i\nu_R}\frac{d \nu_I}{ \nu-k}\, \psi_1(\nu,-2i\nu_I) e_{1_{Xt}}  \tilde\mu_2^+-\frac{1}{\pi}\int_{0}^\infty d\nu_R\int_{-i\infty}^{i\nu_R}\frac{d\nu_I}{\nu-k}\, \psi_2(\nu,-2i\nu_I) e_{1_{Xt}} \tilde \mu_1^+\nonumber\\
&-\frac{1}{\pi}\int_{0}^\infty d\nu_R\int_{-i\infty}^{0}\frac{d\nu_I}{\nu-k}\, \psi_3(\nu,-2i\nu_I) e_{1_{Xt}}  \tilde\mu_2^+-\frac{1}{\pi} \int_{-\infty}^0 d\nu_R\int_{i \nu_R}^{i \infty} \frac{d\nu_I}{\nu-k}\, \chi_1(\nu,-2i\nu_I) e_{1_{Xt}}  \tilde\mu_2^+\nonumber\\
&-\frac{1}{\pi}\int_{0}^\infty d\nu_R\int_{i \nu_R}^{i\infty} \frac{d\nu_I}{\nu-k}\, \chi_2(\nu,-2i\nu_I) e_{1_{Xt}}  \tilde\mu_1^+-\frac{1}{\pi}\int_{-\infty}^0 d\nu_R\int_{i\nu_R}^{0}\frac{d\nu_I}{\nu-k}\, \chi_3(\nu,-2i\nu_I) e_{1_{Xt}} \tilde \mu_2^+\nonumber\\
&-\frac{1}{\pi}\int_{-\infty}^{0}d\nu_R\int_0^\infty\!\! \frac{d\nu_I}{\nu-k}\,f_1^+(\nu_R,\nu_I,T)e_{1_{Xt}}\tilde\mu_2^+-\frac{1}{\pi}\int_0^{\infty}d\nu_R\int_{-\infty}^0 \frac{d\nu_I}{\nu-k}\,f_2^-(\nu_R,\nu_I,T)e_{1_{Xt}}\tilde\mu_1^-.
\end{align}}
Since 
\begin{equation}\label{ordermu}
\mu=1+\mathcal O\left(\frac{1}{k}\right),\quad |k|\rightarrow\infty,
\end{equation}
the first Lax equation \eqref{lax1} implies that
\begin{equation}\label{limit}
 q(x,y,t)=-2i\lim_{|k|\rightarrow \infty}\bigg(k \mu_x(x,y,t,k_R,k_I)\bigg),
\end{equation}
therefore, using equation \eqref{pompi4} and applying l' H\^opital's rule to the above limit, we find equation \eqref{qprop}.\QED
\end{PROOF}

\section{The linear limit}

\begin{figure}[ht]
\begin{center}
\resizebox{7cm}{!}{\input{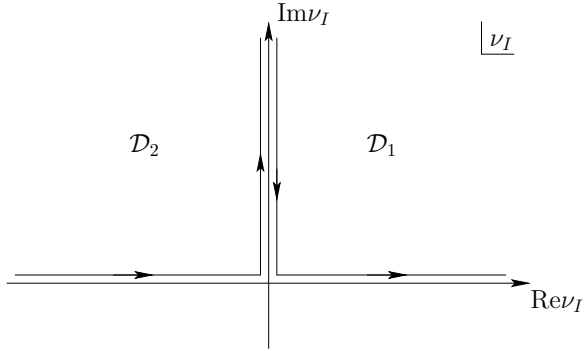}}
\end{center}
\caption{The regions $\mathcal D_1$ and $\mathcal D_2$ on the complex $\nu_I$-plane.}
\label{d1d2}
\end{figure}

In the linear limit
\begin{equation}\label{qu}
q(x,y,t)=\varepsilon u(x,y,t)+\mathcal O(\varepsilon^2), \quad \varepsilon\rightarrow 0,
\end{equation}
the $\mathcal O(\varepsilon)$ term of the KPI equation \eqref{kpiintro} yields the linear PDE
\begin{equation}\label{kpilin}
u_t+u_{xxx}-3 \partial_x^{-1} u_{yy}=0.
\end{equation}
The initial-boundary value problem for this equation can be solved via a spectral analysis similarly to its nonlinear analogue, however, since the PDE is linear we will employ the method of \cite{F1997}. Differentiating equation \eqref{kpilin} with respect to $x$, we have:
\begin{equation}\label{kpilin2}
u_{xt}+u_{xxxx}-3u_{yy}=0,
\end{equation}
with formal adjoint equation for some function $\tilde u$:
\begin{equation}\label{kpilin2adj}
\tilde u_{xt}+\tilde u_{xxxx}-3\tilde u_{yy}=0.
\end{equation}
Multiplying equation \eqref{kpilin2} by $\tilde u_x$, equation \eqref{kpilin2adj} by $u_x$ and adding the two resulting equations, we find:
\begin{equation}\label{divlin}
\left(u_x\tilde u_{x}\right)_t+\left(\tilde u_{x}u_{xxx}+u_x\tilde u_{xxx}-\tilde u_{xx}u_{xx}+3\tilde u_yu_y\right)_x-3\left(\tilde u_{x}u_y+u_x \tilde u_y\right)_y=0.
\end{equation}
A family of solutions to the formal adjoint equation is given by
\begin{equation}\label{fam}
\tilde u(x,y,t)=e^{-ik_1x-ik_2y-i\big(k_1^3+3 \frac{k^2_2}{k_1}\big)t},
\end{equation}
hence
\begin{equation} \label{divlin2}
\begin{split}
& \left(-ik_1 e^{-ik_1x-ik_2y-\big(k_1^3+3 \frac{k^2_2}{k_1}\big)t}u_x\right)_t+3\left(ik_1 u_y+ik_2u_x\right)_y\\
&+\Big(e^{-ik_1x-ik_2y-i\big(k_1^3+3 \frac{k^2_2}{k_1}\big)t}\left(-ik_1 u_{xxx}+i k_1^3u_x+k_1^2 u_{xx}-3i k_2 u_y\right)\Big)_x=0.
\end{split}
\end{equation}
Integrating this expression with respect to $x$ from $-\infty$ to $\infty$ and with respect to $y$ from $0$ to $\infty$, we have:
\begin{equation} \label{divlinint}
\left(k_1^2 e^{-i\big(k_1^3+3 \frac{k^2_2}{k_1}\big)t}\hat u(k_1,k_2,t)\right)_t=3k_1\int_{-\infty}^{\infty}dx\, e^{-ik_1 x-i\big(k_1^3+3 \frac{k^2_2}{k_1}\big)t}\Big[k_2 v(x,t)+i w(x,t)\Big],
\end{equation}
where 
\begin{equation}\label{uhat}
\hat u(k_1,k_2,t):=\int_{-\infty}^{\infty}dx\int_0^\infty dy\, e^{-ik_1x-ik_2y}u(x,y,t)
\end{equation}
and the functions $v$ and $w$ are defined by
\begin{equation}\label{vw}
 v(x,t)=u(x,0,t),\quad w(x,t)=u_y(x,0,t).
\end{equation}
Integrating with respect to $t$ we find the relevant global relation,
\begin{equation}\label{grlin}
\begin{split}
e^{-i\big(k_1^3+3 \frac{k^2_2}{k_1}\big)t}\hat u(k_1,k_2,t)=\hat u_0(k_1,k_2)+\int_0^t d\tau\int_{-\infty}^{\infty}dx\, e^{-ik_1 x-i\big(k_1^3+3 \frac{k^2_2}{k_1}\big)\tau}\, 3\left[\frac{k_2}{k_1}\, v(x,\tau)\right.\\
+\left.\frac{i}{k_1}\, w(x,\tau)\right],\quad k_1\in \mathbb R,\ \mathrm{Im}k_2\leq0.
\end{split}
\end{equation}
Hence, using \eqref{uhat} we obtain:
{\small\begin{align}
&u(x,y,t)=\frac{1}{(2\pi)^2}\int_{-\infty}^{\infty}dk_1\int_{-\infty}^{\infty}dk_2\, e^{ik_1x-ik_2y+i(k_1^3+3\frac{k^2_2}{k_1})t}\hat u_0(k_1,k_2)\nonumber\\
&+\frac{1}{(2\pi)^2}\int_{-\infty}^{\infty}dk_1\int_{-\infty}^{\infty}dk_2\, e^{ik_1x-ik_2y+i(k_1^3+3\frac{k^2_2}{k_1})t}\int_0^t d\tau\int_{-\infty}^{\infty}dx\, e^{-ik_1 x-i(k_1^3+3 \frac{k^2_2}{k_1})\tau}\,\frac{3k_2}{k_1}\, v(x,\tau)\nonumber\\
&+\frac{1}{(2\pi)^2}\int_{-\infty}^{\infty}dk_1\int_{-\infty}^{\infty}dk_2\, e^{ik_1x-ik_2y+i(k_1^3+3\frac{k^2_2}{k_1})t}\int_0^t d\tau\int_{-\infty}^{\infty}dx\, e^{-ik_1 x-i(k_1^3+3 \frac{k^2_2}{k_1})\tau}\,\frac{3i}{k_1}\, w(x,\tau).\label{usol1}
\end{align}}
Furthermore, after changing variables 
\begin{equation}
 k_1=-2\nu_R,\quad k_2=4\nu_R\tilde \nu_I
\end{equation}
and rearranging the integrals, we obtain the following expression:
{\small\begin{align}
&u(x,y,t)=\frac{2}{\pi^2}\mbox{\large $\Big($}\int_{-\infty}^{0}d\nu_R\int_{-\infty}^{\infty}d\tilde \nu_I-\int_{0}^{\infty}d\nu_R\int_{-\infty}^{\infty}d\tilde \nu_I\mbox{\large $\Big)$}\, e^{-2i\nu_R x+4i\nu_R\tilde \nu_Iy-i(8\nu_R^3+24\nu_R\tilde \nu_I^2)t}\hat u_0(-2\nu_R,4\nu_R\tilde \nu_I)\nonumber\\
&-\frac{2}{\pi^2}\int_0^{\infty}\!\!d\nu_R\int_{-\infty}^\infty \!\!d\tilde \nu_I\int_0^t \!\!d\tau\int_{-\infty}^{\infty}\!\!d\xi\, e^{-2i\nu_R (\xi-x)+4i\nu_R\tilde \nu_I y+i(8\nu_R^3+24\nu_R\tilde \nu_I^2)(\tau-t)}\,3\nu_R\mbox{\large $\Big[$}2\tilde \nu_I\, v(\xi,\tau)+\frac{i}{2\nu_R}\, w(\xi,\tau)\mbox{\large $\Big]$}\nonumber\\
&-\frac{2}{\pi^2}\int_{-\infty}^{0}\!\!d\nu_R\int_{-\infty}^\infty\!\! d\tilde \nu_I\int_0^t \!\! d\tau\int_{-\infty}^{\infty}\!\! d\xi\, e^{-2i\nu_R (\xi-x)+4i\nu_R\tilde \nu_I y+i(8\nu_R^3+24\nu_R\tilde \nu_I^2)(\tau-t)}\,3\nu_R\mbox{\large $\Big[$}2\tilde \nu_I\, v(\xi,\tau)+\frac{i}{2\nu_R}\, w(\xi,\tau)\mbox{\large $\Big]$}.\label{usol3}
\end{align}}
Finally, letting $\tilde \nu_I=-i \nu_I$ and recalling the definition \eqref{omega} of $\omega$, we obtain:
{\small\begin{align}
&u(x,y,t)=\frac{2i}{\pi^2}\mbox{\large $\Big($}\int_{-\infty}^{0}d\nu_R\int_{-i\infty}^{i\infty}d \nu_I-\int_{0}^{\infty}d\nu_R\int_{-i\infty}^{i\infty}d \nu_I\mbox{\large $\Big)$}\, e^{-2i\nu_R x+4\nu_R \nu_Iy+\omega(\nu,-2\nu_R)t}\hat u_0(-2\nu_R,-4i\nu_R \nu_I)\nonumber\\
&-\frac{2}{\pi^2}\int_0^{\infty}\!\!d\nu_R\int_{-i\infty}^{i\infty} \!\!d \nu_I\int_0^t \!\!d\tau\int_{-\infty}^{\infty}\!\!d\xi\, e^{-2i\nu_R (\xi-x)+4\nu_R \nu_I y-\omega(\nu,-2\nu_R)(\tau-t)}\,3\nu_R\mbox{\large $\Big[$}2 \nu_I\, v(\xi,\tau)-\frac{1}{2\nu_R}\, w(\xi,\tau)\mbox{\large $\Big]$}\nonumber\\
&-\frac{2}{\pi^2}\int_{-\infty}^{0}\!\!d\nu_R\int_{-i\infty}^{i\infty}\!\! d \nu_I\int_0^t \!\! d\tau\int_{-\infty}^{\infty}\!\! d\xi\, e^{-2i\nu_R (\xi-x)+4\nu_R \nu_I y-\omega(\nu,-2\nu_R)(\tau-t)}\,3\nu_R\mbox{\large $\Big[$}2 \nu_I\, v(\xi,\tau)-\frac{1}{2\nu_R}\, w(\xi,\tau)\mbox{\large $\Big]$}.\label{usol4}
\end{align}}
By employing Cauchy's theorem and Jordan's lemma, the $\nu_I$-contours of the second and the third term on the RHS of the above expression can be deformed to the positively oriented boundaries of the first and the second quadrant of the complex $\nu_I$-plane respectively (see figure \ref{d1d2}), i.e.

{\small\begin{align}
&u(x,y,t)=\frac{2i}{\pi^2}\mbox{\large $\Big($}\int_{-\infty}^{0}d\nu_R\int_{-i\infty}^{i\infty}d \nu_I-\int_{0}^{\infty}d\nu_R\int_{-i\infty}^{i\infty}d \nu_I\mbox{\large $\Big)$}\, e^{-2i\nu_R x+4\nu_R \nu_Iy+\omega(\nu,-2\nu_R)t}\hat u_0(-2\nu_R,-4i\nu_R \nu_I)\nonumber\\
&-\frac{2}{\pi^2}\int_0^{\infty}\!\!d\nu_R\int_{\partial\mathcal D_2} \!\!d \nu_I\int_0^t \!\!d\tau\int_{-\infty}^{\infty}\!\!d\xi\, e^{-2i\nu_R (\xi-x)+4\nu_R \nu_I y-\omega(\nu,-2\nu_R)(\tau-t)}\,3\nu_R\mbox{\large $\Big[$}2 \nu_I\, v(\xi,\tau)-\frac{1}{2\nu_R}\, w(\xi,\tau)\mbox{\large $\Big]$}\nonumber\\
&-\frac{2}{\pi^2}\int_{-\infty}^{0}\!\!d\nu_R\int_{\partial \mathcal D_1}\!\! d \nu_I\int_0^t \!\! d\tau\int_{-\infty}^{\infty}\!\! d\xi\, e^{-2i\nu_R (\xi-x)+4\nu_R \nu_I y-\omega(\nu,-2\nu_R)(\tau-t)}\,3\nu_R\mbox{\large $\Big[$}2 \nu_I\, v(\xi,\tau)-\frac{1}{2\nu_R}\, w(\xi,\tau)\mbox{\large $\Big]$}.\label{usol5}
\end{align}}
Moreover, Jordan's lemma implies that
{\small\begin{align}
&\int_0^{\infty}\!\!d\nu_R\int_{\partial\mathcal D_2} \!\!d \nu_I\int_t^T \!\!d\tau\int_{-\infty}^{\infty}\!\!d\xi\, e^{-2i\nu_R (\xi-x)+4\nu_R \nu_I y-\omega(\nu,-2\nu_R)(\tau-t)}\, (\cdot)=\nonumber\\
&=\int_{-\infty}^{0}\!\!d\nu_R\int_{\partial\mathcal D_1}\!\! d \nu_I\int_t^T \!\! d\tau\int_{-\infty}^{\infty}\!\! d\xi\, e^{-2i\nu_R (\xi-x)+4\nu_R \nu_I y-\omega(\nu,-2\nu_R)(\tau-t)}\,(\cdot)=0,
\end{align}}
thus the upper limit of the $\tau$-integrals can be replaced by $T$, hence  
{\small\begin{align}
&u(x,y,t)=\frac{2i}{\pi^2}\mbox{\large $\Big($}\int_{-\infty}^{0}d\nu_R\int_{-i\infty}^{i\infty}d \nu_I-\int_{0}^{\infty}d\nu_R\int_{-i\infty}^{i\infty}d \nu_I\mbox{\large $\Big)$}\, e^{-2i\nu_R x+4\nu_R \nu_Iy+\omega(\nu,-2\nu_R)t}\hat u_0(-2\nu_R,-4i\nu_R \nu_I)\nonumber\\
&-\frac{2}{\pi^2}\int_0^{\infty}\!\!d\nu_R\int_{\partial\mathcal D_2} \!\!d \nu_I\int_0^T \!\!d\tau\int_{-\infty}^{\infty}\!\!d\xi\, e^{-2i\nu_R (\xi-x)+4\nu_R \nu_I y-\omega(\nu,-2\nu_R)(\tau-t)}\,3\nu_R\mbox{\large $\Big[$}2 \nu_I\, v(\xi,\tau)-\frac{1}{2\nu_R}\, w(\xi,\tau)\mbox{\large $\Big]$}\nonumber\\
&-\frac{2}{\pi^2}\int_{-\infty}^{0}\!\!d\nu_R\int_{\partial\mathcal D_1}\!\! d \nu_I\int_0^T \!\! d\tau\int_{-\infty}^{\infty}\!\! d\xi\, e^{-2i\nu_R (\xi-x)+4\nu_R \nu_I y-\omega(\nu,-2\nu_R)(\tau-t)}\,3\nu_R\mbox{\large $\Big[$}2 \nu_I\, v(\xi,\tau)-\frac{1}{2\nu_R}\, w(\xi,\tau)\mbox{\large $\Big]$}.\label{usol6}
\end{align}}

On the other hand, taking the linear limit directly on the representation \eqref{pompi4} obtained for $\mu$ via the spectral analysis of the Lax pair of the KPI equation, i.e. letting
\begin{equation}\label{linlim}
\mu=1+\varepsilon M+\mathcal O(\varepsilon^2), \quad r_j^\pm=\varepsilon R+\mathcal O(\varepsilon^2), \quad p_j^\pm=\varepsilon P+\mathcal O(\varepsilon^2),\quad j=1,2, \quad \varepsilon\rightarrow 0,
\end{equation}
implies via equations \eqref{chi} and \eqref{psi}:
\begin{equation}\label{linlim2}
\chi_j=\varepsilon X_j+\mathcal O(\varepsilon^2), \quad \psi_j=\varepsilon \Psi_j+\mathcal O(\varepsilon^2), \quad j=1,2,3,
\end{equation}
for the functions $X_j$ and $\Psi_j$ defined by: 
\begin{align}
&X_1(k_R,\lambda)=\Psi_1(k_R,\lambda)=-R(k_R,-2k_R-\lambda),\nonumber\\
& X_2(k_R,\lambda)=\Psi_2(k_R,\lambda)=P(k_R,-2k_R-\lambda) +R(k_R,-2k_R-\lambda),\nonumber\\
& X_3(k_R,\lambda)=\Psi_3(k_R,\lambda)=-P(k_R,-2k_R-\lambda), \label{X}
\end{align}
where
\begin{align}
R(k,l)&=\frac{i}{2\pi}\int_{-\infty}^{\infty}d\xi\int_{0}^{\infty}d\eta\, e^{-il\xi-il(l+2k)\eta}u_0(\xi,\eta)\label{R},\\
P(k,l)&=\frac{1}{2\pi}\int_{-\infty}^{\infty}d\xi\int_{0}^{T}d\tau\, e^{-il\xi+\omega(k,l)\tau}\, 3\Big[-i(l+2k)v(\xi,\tau)-i\partial_x^{-1}w(\xi,\tau)\Big].\label{P}
\end{align}
Hence, the $\mathcal O(\varepsilon)$ term of equation \eqref{pompi4} yields:
{\small\begin{align}
M&=\frac{1}{\pi}\int_{-\infty}^{0}d\nu_R\int_{-i\infty}^{i\nu_R}\frac{d\nu_I}{\nu-k}\, e_{1_{Xt}} R(\nu,-2\nu_R)\nonumber\\
&-\frac{1}{\pi}\int_{0}^{\infty}d\nu_R\int_{-i\infty}^{i\nu_R}\frac{d\nu_I}{\nu-k}\, e_{1_{Xt}} \mbox{\large $\Big[$} P(\nu,-2\nu_R)+R(\nu,-2\nu_R)\mbox{\large $\Big]$}\nonumber\\
&+\frac{1}{\pi}\int_{0}^{\infty}d\nu_R\int_{-i\infty}^{0}\frac{d\nu_I}{\nu-k}\, e_{1_{Xt}} P(\nu,-2\nu_R)+\frac{1}{\pi}\int_{-\infty}^{0}d\nu_R\int_{i\nu_R}^{i\infty}\frac{d\nu_I}{\nu-k}\, e_{1_{Xt}} R(\nu,-2\nu_R)\nonumber\\
&-\frac{1}{\pi}\int_{0}^{\infty}d\nu_R\int_{i\nu_R}^{i\infty}\frac{d\nu_I}{\nu-k}\, e_{1_{Xt}} \mbox{\large $\Big[$}P(\nu,-2\nu_R)+R(\nu,-2\nu_R)\mbox{\large $\Big]$}\nonumber\\
&+\frac{1}{\pi}\int_{-\infty}^{0}d\nu_R\int_{0}^{i\infty}\frac{d\nu_I}{\nu-k}\, e_{1_{Xt}} P(\nu,-2\nu_R)-\frac{1}{\pi}\int_{-\infty}^{0}d\nu_R\int_{0}^{\infty}\frac{d\nu_I}{\nu-k}\, e_{1_{Xt}} P(\nu,-2\nu_R)\nonumber\\
&-\frac{1}{\pi}\int_{0}^{\infty}d\nu_R\int_{-\infty}^{0}\frac{d\nu_I}{\nu-k}\, e_{1_{Xt}} P(\nu,-2\nu_R).\label{M}
\end{align}}
Thus, taking the limit \eqref{limit} we find:
{\small\begin{align}
u&=\frac{2i}{\pi}\int_{-\infty}^{0}d\nu_R\int_{-i\infty}^{i\infty}d\nu_I\, e_{1_{Xt}} (-2i\nu_R)R(\nu,-2\nu_R)\nonumber\\
&-\frac{2i}{\pi}\int_{0}^{\infty}d\nu_R\int_{-i\infty}^{i\infty}d\nu_I\, e_{1_{Xt}} (-2i\nu_R)\mbox{\large $\Big[$} P(\nu,-2\nu_R)+R(\nu,-2\nu_R)\mbox{\large $\Big]$}\nonumber\\
&-\frac{2i}{\pi}\int_{0}^{\infty}d\nu_R\int_{-\infty}^{0}d\nu_I\, e_{1_{Xt}}(-2i\nu_R) P(\nu,-2\nu_R)-\frac{2i}{\pi}\int_{0}^{\infty}d\nu_R\int_{0}^{-i\infty}d\nu_I\, e_{1_{Xt}}(-2i\nu_R) P(\nu,-2\nu_R)\nonumber\\
&+\frac{2i}{\pi}\int_{-\infty}^{0}d\nu_R\int_{0}^{i\infty}d\nu_I\, e_{1_{Xt}}(-2i\nu_R) P(\nu,-2\nu_R)-\frac{2i}{\pi}\int_{-\infty}^0 d\nu_R\int_{0}^{\infty}d\nu_I\, e_{1_{Xt}} (-2i\nu_R)P(\nu,-2\nu_R),\label{u1}
\end{align}}
or, using the definition \eqref{eXt} and integrating once by parts with respect to $x$, 
{\small\begin{align}
u&=\frac{2i}{\pi^2}\mbox{\large $\Big($}\int_{-\infty}^{0}\!\!d\nu_R\int_{-i\infty}^{i\infty}d \nu_I-\int_{0}^{\infty}\!\!d\nu_R\int_{-i\infty}^{i\infty}d \nu_I\mbox{\large $\Big)$}\!\!\int_{-\infty}^{\infty}\!\!d\xi\int_{0}^{\infty}\!\!d\eta\, e^{2i\nu_R (\xi-x)-4\nu_R \nu_I(\eta-y)+\omega(\nu,-2\nu_R)t}\nu_R\, u_0(\xi,\eta)\nonumber\\
&-\frac{2}{\pi^2}\int_0^{\infty}\!\!d\nu_R\int_{\partial\mathcal D_2}\!\!d \nu_I \!\!\int_0^t \!\!d\tau\!\!\int_{-\infty}^{\infty}\!\!d\xi\, e^{-2i\nu_R (\xi-x)+4\nu_R \nu_I y-\omega(\nu,-2\nu_R)(\tau-t)}\,3\nu_R\mbox{\large $\Big[$}2 \nu_I\, v(\xi,\tau)-\frac{1}{2\nu_R}\, w(\xi,\tau)\mbox{\large $\Big]$}\nonumber\\
&-\frac{2}{\pi^2}\int_{-\infty}^{0}\!\!d\nu_R\int_{\partial\mathcal {D_1}}\!\!d \nu_I\!\!\int_0^t \!\! d\tau\int_{-\infty}^{\infty}\!\! d\xi\, e^{-2i\nu_R (\xi-x)+4\nu_R \nu_I y-\omega(\nu,-2\nu_R)(\tau-t)}\,3\nu_R\mbox{\large $\Big[$}2 \nu_I\, v(\xi,\tau)-\frac{1}{2\nu_R}\, w(\xi,\tau)\mbox{\large $\Big]$}.\label{u2}
\end{align}}
Recalling the definition \eqref{uhat}, this equation is precisely equation \eqref{usol6}. 

\section{Conclusion}

The integral representation \eqref{qprop} for the solution $q(x,y,t)$ of the initial-boundary value problem for the KPI equation \eqref{kpiintro} on the half-plane has been constructed, following the approach to the initial-boundary value problems for the DS and the KPII equations presented in \cite{FDS2009} and \cite{MF2011}. In particular, the spectral analysis of the associated Lax pair of the PDE is achieved by using several tools of the theory of complex variables, including Green's theorem and the so-called Pompeiu's formula, in the context of a d-bar formalism.

In addition, the important identity \eqref{griprop} is obtained which relates certain transforms of the initial condition $q_0(x,y)$ and the boundary values $g(x,t)$ and $h(x,t)$, called the \textit{spectral functions}, to the corresponding transform of the solution $q(x,y,t)$, and which is valid in appropriate regions of the complex $k$-plane. This identity, called the \textit{global relation}, can in principle be employed after appropriate transformations in the spectral variable  $k$ in order to eliminate the spectral function of the \textit{unknown} boundary value (either $g$ or $h$) from the representation \eqref{qprop}. This is a rather complicated task, since both \eqref{griprop} and \eqref{qprop} involve the function $\mu(x,y,t,k_R,k_I)$, which is given by equation \eqref{pompi4} as the solution of certain linear Volterra integral equations. However, it has been shown (see for example \cite{FSG2004}, \cite{FLenells2010b}) that for a particular class of boundary conditions, called \textit{linearisable}, one can obtain an expression for $q$ depending only on spectral functions of the \textit{given} initial and boundary values. The case of linearisable boundary conditions for the KP equations is currently being investigated.

\section*{Acknowledgements}
\noindent D. Mantzavinos was supported by an EPSRC Doctoral Training Grant.


\begin{thebibliography}{111}
\bibitem{GGKM} C. S. Gardner, J. M. Greene, M. D. Kruskal, and R. M. Miura, Method for Solving the Korteweg-de Vries Equation, Phys. Rev. Lett. \textbf{19}, 1095-1097 (1967).
\bibitem{FAbl1983} A. S. Fokas and M. J. Ablowitz, On the Inverse Scattering of the Time Dependent Equation and the Associated KPI Equation, Stud. Appl. Math. \textbf{69}, 211-228 (1983).
\bibitem{AblBarF1983} M. J. Ablowitz, D. BarYaacov and A. S. Fokas, On the Inverse Scattering Transform for the Kadomtsev-Petvisvhili Equation, Stud. Appl. Math. \textbf{69}, 135-143 (1983). 
\bibitem{BLP1989} M. Boiti, J.J.-P. Leon, F. Pempinelli, Spectral Transform and Orthogonality Relations for the Kadomtsev-Petviashvili Equation, Physics Letters A, \textbf{141}, 96-100 (1989).
\bibitem{FSung1992} A. S. Fokas and L. Y. Sung, On the Solvability of the N-Wave, the Davey-Stewartson and the Kadomtsev-Petviashvili Equation, Inverse Problems \textbf{8}, 673-708 (1992).
\bibitem{BC1984} R. Beals and R.R. Coifman, Linear Spectral Problem, Nonlinear Equations and the $\bar \partial$ Method, Inverse Problem {\bf 5}, 87-130 (1989).
\bibitem{Fkp2009}  A. S. Fokas, The Kadomtsev-Petviashvili Equation Revisited and Integrability in 4+2 and 3+1, Stud. Appl. Math. \textbf{122}, 347-359 (2009).
\bibitem{F1997} A. S. Fokas, A Unified Transform Method for Solving Linear and Certain Nonlinear PDEs, Proc. R. Soc. Lond. A \textbf{453}, 1411-1443 (1997).
\bibitem{F2002a} A. S. Fokas, Integrable Nonlinear Evolution Equations on the Half-Line, Comm. Math. Phys. \textbf{230}, 1-39 (2002).
\bibitem{Fbook} A. S. Fokas, \emph{A Unified Approach to Boundary Value Problems}, SIAM (2008).
\bibitem{F2002b} A. S. Fokas, A New Transform Method for Evolution PDEs, IMA J. Appl. Math. \textbf{67}, 1-32 (2002).
\bibitem{FTreharne2008} P. A. Treharne, A. S. Fokas, The Generalized Dirichlet to Neumann map for the KdV Equation on the Half-Line, J. Nonlinear Science \textbf{18}, 191-217 (2008).
\bibitem{FDS2009} A. S. Fokas, The Davey-Stewartson Equation on the Half-Plane, Comm. Math. Phys. \textbf{289}, 957-993 (2009).
\bibitem{MF2011} D. Mantzavinos and A. S. Fokas, The Kadomtsev-Petviashvili II equation on the half-plane, Physica D \textbf{240}, 477-511 (2011).
\bibitem{FSG2004} A. S. Fokas, Linearizable Initial-Boundary Value Problems for the sine-Gordon Equation on the Half-Line, Nonlinearity \textbf{17}, 1521 (2004).
\bibitem{FLenells2010b} A. S. Fokas and J. Lenells, Linearizable Boundary Value Problems for the Elliptic sine-Gordon and the Elliptic Ernst Equations (submitted).


\end{thebibliography}
\end{document}